\documentclass[journal, twoside]{IEEEtran}

\usepackage{setspace}
\usepackage{relsize}
\usepackage{url}

\usepackage{amsfonts}
\usepackage[cmex10]{amsmath} % Needs to be defined before amsthm package
\usepackage{amssymb}

\usepackage{amsthm}	% Needs to be defined after amsmath package
\usepackage{mathrsfs}	% Ralph Smith’s Formal Script Font (choose between mathrsfs or euscript for using \mathscr{})
\usepackage[pdftex]{graphicx}	% Option pdftex added as recommended for epstopdf package below
\usepackage{epstopdf}	% Added to support eps files with PDFLaTex
\usepackage[breaklinks=true,linkcolor=blue,anchorcolor=black,citecolor=blue,urlcolor=blue]{hyperref}
\usepackage{color}
\usepackage{authblk}
\usepackage{cite}
\usepackage{multirow}
\usepackage[font=footnotesize]{caption}
\usepackage[font=footnotesize,subrefformat=parens]{subcaption}
\usepackage{hyperref}
\usepackage{array}									% Generic table/tabular package minimally required for IEEE journals

\theoremstyle{plain} % {plain}, definition, remark (bold title, italic text)
\newtheorem{theorem}{Theorem}
\newtheorem{proposition}{Proposition}
\newtheorem{corollary}{Corollary}

\theoremstyle{definition} % plain, {definition}, remark (bold title, regular text)
\newtheorem{lemma}{Lemma}
\newtheorem{definition}{Definition}
\newtheorem{problem}{Problem}
\newtheorem{property}{Property}
\newtheorem{conjecture}{Conjecture}
\newtheorem{example}{Example}
\newtheorem{assumption}{Assumption}

\theoremstyle{remark} % plain, definition, {remark} (title in italics, regular text)
\newtheorem{remark}{Remark}
\newtheorem{fact}{Fact}
\newtheorem{claim}{Claim}

\newcommand{\theoremAlt}[2][]{		\begin{theorem}{{\textit{#1}}} 	{#2} \end{theorem}}
\newcommand{\corollaryAlt}[2][]{		\begin{corollary}{\textit{{#1}}}	{#2} \end{corollary}}

\newcommand{\assumptionAlt}[2][]{	\begin{assumption}{{#1}}		{#2} \end{assumption}}

\newcommand{\propositionAlt}[2][]{	\begin{proposition}{{#1}}		{#2} \end{proposition}}

\def\nn{\nonumber}
\def\beq{\begin{equation}}
\def\eeq{\end{equation}}
\def\beqs{\begin{equation*}}
\def\eeqs{\end{equation*}}
\def\balign{\begin{align}}
\def\ealign{\end{align}}
\def\bea{\begin{eqnarray}}
\def\eea{\end{eqnarray}}
\def\ba{\begin{array}}
\def\ea{\end{array}}
\def\bcenter{\begin{center}}
\def\ecenter{\end{center}}
\def\bitem{\begin{itemize}}
\def\eitem{\end{itemize}}
\def\benum{\begin{enumerate}}
\def\eenum{\end{enumerate}}
\def\bdoc{\begin{document}}
\def\edoc{\end{document}}
\def\defeq{{\stackrel{\Delta}{=}}}
\def\nn{\nonumber}

\newcommand{\eg}{\mbox{\textit{e.g.}, \/}}
\newcommand{\etal}{\mbox{\textit{et al.} \/}}
\newcommand{\viz}{\mbox{\textit{\iet viz.},\ \/}}
\newcommand{\ie}{\mbox{\textit{i.e.},\ \/}}
\newcommand{\etc}{\mbox{\textit{etc.},\ \/}}
\newcommand{\cf}{\mbox{\textit{cf.},\ \/}}
\newcommand{\vs}{\mbox{\textit{vs.}\ \/}}

\def\red{\textcolor{red}}
\def\blue{\textcolor{blue}}
\def\black{\textcolor{black}}
\def\yellow{\textcolor{yellow}}
\def\green{\textcolor{green}}
\def\tcbg{\textcolor{bgrd}}
\def\magenta{\textcolor{magenta}}
\def\white{\textcolor{white}}
\def\gray{\textcolor{gray}}
\newcommand\Quote[1]{\lq\textsl{#1}\rq}
\newcommand\fr[2]{{\textstyle\frac{#1}{#2}}}

\newcommand{\Ambb}{\mathbb{A}}
\newcommand{\Bmbb}{\mathbb{B}}
\newcommand{\Cmbb}{\mathbb{C}}
\newcommand{\Dmbb}{\mathbb{D}}
\newcommand{\Embb}{\mathbb{E}}
\newcommand{\Fmbb}{\mathbb{F}}
\newcommand{\Gmbb}{\mathbb{G}}
\newcommand{\Hmbb}{\mathbb{H}}
\newcommand{\Imbb}{\mathbb{I}}
\newcommand{\Jmbb}{\mathbb{J}}
\newcommand{\Kmbb}{\mathbb{K}}
\newcommand{\Lmbb}{\mathbb{L}}
\newcommand{\Mmbb}{\mathbb{M}}
\newcommand{\Nmbb}{\mathbb{N}}
\newcommand{\Ombb}{\mathbb{O}}
\newcommand{\Pmbb}{\mathbb{P}}
\newcommand{\Qmbb}{\mathbb{Q}}
\newcommand{\Rmbb}{\mathbb{R}}
\newcommand{\Smbb}{\mathbb{S}}
\newcommand{\Tmbb}{\mathbb{T}}
\newcommand{\Umbb}{\mathbb{U}}
\newcommand{\Vmbb}{\mathbb{V}}
\newcommand{\Wmbb}{\mathbb{W}}
\newcommand{\Xmbb}{\mathbb{X}}
\newcommand{\Ymbb}{\mathbb{Y}}
\newcommand{\Zmbb}{\mathbb{Z}}

\newcommand{\Amsc}{\mathcal{A}}
\newcommand{\Bmsc}{\mathcal{B}}
\newcommand{\Cmsc}{\mathcal{C}}
\newcommand{\Dmsc}{\mathcal{D}}
\newcommand{\Emsc}{\mathcal{E}}
\newcommand{\Fmsc}{\mathcal{F}}
\newcommand{\Gmsc}{\mathcal{G}}
\newcommand{\Hmsc}{\mathcal{H}}
\newcommand{\Imsc}{\mathcal{I}}
\newcommand{\Jmsc}{\mathcal{J}}
\newcommand{\Kmsc}{\mathcal{K}}
\newcommand{\Lmsc}{\mathcal{L}}
\newcommand{\Mmsc}{\mathcal{M}}
\newcommand{\Nmsc}{\mathcal{N}}
\newcommand{\Omsc}{\mathcal{O}}
\newcommand{\Pmsc}{\mathcal{P}}
\newcommand{\Qmsc}{\mathcal{Q}}
\newcommand{\Rmsc}{\mathcal{R}}
\newcommand{\Smsc}{\mathcal{S}}
\newcommand{\Tmsc}{\mathcal{T}}
\newcommand{\Umsc}{\mathcal{U}}
\newcommand{\Vmsc}{\mathcal{V}}
\newcommand{\Wmsc}{\mathcal{W}}
\newcommand{\Xmsc}{\mathcal{X}}
\newcommand{\Ymsc}{\mathcal{Y}}
\newcommand{\Zmsc}{\mathcal{Z}}

\newcommand{\Ascr}{\mathscr{A}}
\newcommand{\Bscr}{\mathscr{B}}
\newcommand{\Cscr}{\mathscr{C}}
\newcommand{\Dscr}{\mathscr{D}}
\newcommand{\Escr}{\mathscr{E}}
\newcommand{\Fscr}{\mathscr{F}}
\newcommand{\Gscr}{\mathscr{G}}
\newcommand{\Hscr}{\mathscr{H}}
\newcommand{\Iscr}{\mathscr{I}}
\newcommand{\Jscr}{\mathscr{J}}
\newcommand{\Kscr}{\mathscr{K}}
\newcommand{\Lscr}{\mathscr{L}}
\newcommand{\Mscr}{\mathscr{M}}
\newcommand{\Nscr}{\mathscr{N}}
\newcommand{\Oscr}{\mathscr{O}}
\newcommand{\Pscr}{\mathscr{P}}
\newcommand{\Qscr}{\mathscr{Q}}
\newcommand{\Rscr}{\mathscr{R}}
\newcommand{\Sscr}{\mathscr{S}}
\newcommand{\Tscr}{\mathscr{T}}
\newcommand{\Uscr}{\mathscr{U}}
\newcommand{\Vscr}{\mathscr{V}}
\newcommand{\Wscr}{\mathscr{W}}
\newcommand{\Xscr}{\mathscr{X}}
\newcommand{\Yscr}{\mathscr{Y}}
\newcommand{\Zscr}{\mathscr{Z}}

\def\Atiny{\mbox{\tiny{A}}}
\def\Btiny{\mbox{\tiny{B}}}
\def\Ctiny{\mbox{\tiny{C}}}
\def\Dtiny{\mbox{\tiny{D}}}
\def\Etiny{\mbox{\tiny{E}}}
\def\Ctiny{\mbox{\tiny{C}}}
\def\Ftiny{\mbox{\tiny{F}}}
\def\Gtiny{\mbox{\tiny{G}}}
\def\Htiny{\mbox{\tiny{H}}}
\def\Itiny{\mbox{\tiny{I}}}
\def\Jtiny{\mbox{\tiny{J}}}
\def\Ktiny{\mbox{\tiny{K}}}
\def\Ltiny{\mbox{\tiny{L}}}
\def\Mtiny{\mbox{\tiny{M}}}
\def\Ntiny{\mbox{\tiny{N}}}
\def\Otiny{\mbox{\tiny{O}}}
\def\Ptiny{\mbox{\tiny{P}}}
\def\Qtiny{\mbox{\tiny{Q}}}
\def\Rtiny{\mbox{\tiny{R}}}
\def\Stiny{\mbox{\tiny{S}}}
\def\Ttiny{\mbox{\tiny{T}}}
\def\Utiny{\mbox{\tiny{U}}}
\def\Vtiny{\mbox{\tiny{V}}}
\def\Wtiny{\mbox{\tiny{W}}}
\def\Xtiny{\mbox{\tiny{X}}}
\def\Ytiny{\mbox{\tiny{Y}}}
\def\Ztiny{\mbox{\tiny{Z}}}

\def\alphabf{\hbox{\boldmath$\alpha$\unboldmath}}
\def\betabf{\hbox{\boldmath$\beta$\unboldmath}}
\def\gammabf{\hbox{\boldmath$\gamma$\unboldmath}}
\def\deltabf{\hbox{\boldmath$\delta$\unboldmath}}
\def\epsilonbf{\hbox{\boldmath$\epsilon$\unboldmath}}
\def\zetabf{\hbox{\boldmath$\zeta$\unboldmath}}
\def\etabf{\hbox{\boldmath$\eta$\unboldmath}}
\def\iotabf{\hbox{\boldmath$\iota$\unboldmath}}
\def\kappabf{\hbox{\boldmath$\kappa$\unboldmath}}
\def\lambdabf{\hbox{\boldmath$\lambda$\unboldmath}}
\def\mubf{\hbox{\boldmath$\mu$\unboldmath}}
\def\nubf{\hbox{\boldmath$\nu$\unboldmath}}
\def\xibf{\hbox{\boldmath$\xi$\unboldmath}}
\def\pibf{\hbox{\boldmath$\pi$\unboldmath}}
\def\rhobf{\hbox{\boldmath$\rho$\unboldmath}}
\def\sigmabf{\hbox{\boldmath$\sigma$\unboldmath}}
\def\taubf{\hbox{\boldmath$\tau$\unboldmath}}
\def\upsilonbf{\hbox{\boldmath$\upsilon$\unboldmath}}
\def\phibf{\hbox{\boldmath$\phi$\unboldmath}}
\def\chibf{\hbox{\boldmath$\chi$\unboldmath}}
\def\psibf{\hbox{\boldmath$\psi$\unboldmath}}
\def\omegabf{\hbox{\boldmath$\omega$\unboldmath}}
\def\Sigmabf{\hbox{$\bf \Sigma$}}
\def\Upsilonbf{\hbox{$\bf \Upsilon$}}
\def\Omegabf{\hbox{$\bf \Omega$}}
\def\Deltabf{\hbox{$\bf \Delta$}}
\def\Gammabf{\hbox{$\bf \Gamma$}}
\def\Thetabf{\hbox{$\bf \Theta$}}
\def\Lambdabf{\mbox{$ \bf \Lambda $}}
\def\Xibf{\hbox{\bf$\Xi$}}
\def\Pibf{{\bf \Pi}}
\def\thetabf{{\mbox{\boldmath$\theta$\unboldmath}}}
\def\Upsilonbf{\hbox{\boldmath$\Upsilon$\unboldmath}}
\newcommand{\Phibf}{\mbox{${\bf \Phi}$}}
\newcommand{\Psibf}{\mbox{${\bf \Psi}$}}

\def\abf{{\bf a}}
\def\bbf{{\bf b}}
\def\cbf{{\bf c}}
\def\dbf{{\bf d}}
\def\ebf{{\bf e}}
\def\fbf{{\bf f}}
\def\gbf{{\bf g}}
\def\hbf{{\bf h}}
\def\ibf{{\bf i}}
\def\jbf{{\bf j}}
\def\kbf{{\bf k}}
\def\lbf{{\bf l}}
\def\mbf{{\bf m}}
\def\nbf{{\bf n}}
\def\obf{{\bf o}}
\def\pbf{{\bf p}}
\def\qbf{{\bf q}}
\def\rbf{{\bf r}}
\def\sbf{{\bf s}}
\def\tbf{{\bf t}}
\def\ubf{{\bf u}}
\def\vbf{{\bf v}}
\def\wbf{{\bf w}}
\def\xbf{{\bf x}}
\def\ybf{{\bf y}}
\def\zbf{{\bf z}}
\def\rbf{{\bf r}}
\def\xbf{{\bf x}}
\def\ybf{{\bf y}}
\def\Abf{{\bf A}}
\def\Bbf{{\bf B}}
\def\Cbf{{\bf C}}
\def\Dbf{{\bf D}}
\def\Ebf{{\bf E}}
\def\Fbf{{\bf F}}
\def\Gbf{{\bf G}}
\def\Hbf{{\bf H}}
\def\Ibf{{\bf I}}
\def\Jbf{{\bf J}}
\def\Kbf{{\bf K}}
\def\Lbf{{\bf L}}
\def\Mbf{{\bf M}}
\def\Nbf{{\bf N}}
\def\Obf{{\bf O}}
\def\Pbf{{\bf P}}
\def\Qbf{{\bf Q}}
\def\Rbf{{\bf R}}
\def\Sbf{{\bf S}}
\def\Tbf{{\bf T}}
\def\Ubf{{\bf U}}
\def\Vbf{{\bf V}}
\def\Wbf{{\bf W}}
\def\Xbf{{\bf X}}
\def\Ybf{{\bf Y}}
\def\Zbf{{\bf Z}}
\def\Ac{{\cal A}}
\def\Bc{{\cal B}}
\def\Cc{{\cal C}}
\def\Dc{{\cal D}}
\def\Ec{{\cal E}}
\def\Fc{{\cal F}}
\def\Gc{{\cal G}}
\def\Hc{{\cal H}}
\def\Ic{{\cal I}}
\def\Jc{{\cal J}}
\def\Kc{{\cal K}}
\def\Lc{{\cal L}}
\def\Mc{{\cal M}}
\def\Nc{{\cal N}}
\def\Oc{{\cal O}}
\def\Pc{{\cal P}}
\def\Qc{{\cal Q}}
\def\Rc{{\cal R}}
\def\Sc{{\cal S}}
\def\Tc{{\cal T}}
\def\Uc{{\cal U}}
\def\Vc{{\cal V}}
\def\Wc{{\cal W}}
\def\Xc{{\cal X}}
\def\Yc{{\cal Y}}
\def\Zc{{\cal Z}}

\def\0bf{\mathbf{0}}
\def\1bf{\mathbf{1}}

\def\lambdabf{\boldsymbol{\lambda}}
\def\mubf{\boldsymbol{\mu}}
\def\gammabf{\boldsymbol{\gamma}}
\def\xibf{\boldsymbol{\xi}}
\def\sigmabf{\boldsymbol{\sigma}}
\def\zetabf{\boldsymbol{\zeta}}

\def\Lambdabf{\boldsymbol{\Lambda}}
\def\Mubf{\boldsymbol{\Mu}}
\def\Gammabf{\boldsymbol{\Gamma}}
\def\Xibf{\boldsymbol{\Xi}}
\def\Sigmabf{\boldsymbol{\Sigma}}
\def\Zetabf{\boldsymbol{\Zeta}}

\def\1bf{\mathbf{1}}

\def\eqdef{\overset{\text{{\tiny def}}}{=}}

\DeclareMathOperator*{\cart}{\times}
\DeclareMathOperator*{\cov}{\Cmbb\textup{ov}}
\DeclareMathOperator*{\conv}{\text{conv}}
\DeclareMathOperator*{\diag}{\text{diag}}
\DeclareMathOperator*{\tr}{\textup{Tr}}
\DeclareMathOperator*{\minimize}{\text{minimize}}
\DeclareMathOperator*{\proj}{\text{proj}}
\DeclareMathOperator*{\toeplitz}{\text{toeplitz}}
\DeclareMathOperator*{\corr}{\text{corr}}
\DeclareMathOperator*{\rows}{\text{rows}}

\newcommand{\sfrac}[2]{\mbox{$\frac{#1}{#2}$}}

\let\vec\relax
\DeclareMathOperator*{\vec}{\text{vec}}

\newcommand{\Top}{{\mbox{{\tiny $\top$}}}}

\renewcommand{\b}[1]{\ensuremath{\boldsymbol{\mathrm{#1}}}}

\def\sw{\textup{\textsf{sw}}}
\def\rp{\textup{\textsf{rs}}}
\def\cs{\textup{\textsf{cs}}}
\def\so{\textup{\textsf{so}}}
\def\vp{\textup{\textsf{vp}}}

\def\figwid{4cm} %5.7cm for 3 fig
\newcommand*{\MyLocalPath}{}%

\newcommand{\ifThesis}[2]{\ifdefined\IEEEPARstart%
#2%
\else%
#1%
\fi}

\renewcommand{\qedsymbol}{\blacksquare}
\newcommand{\qedadhoc}{\tag*{$\blacksquare$}}

\providecommand{\boldsymbol}[1]{\mbox{\boldmath $#1$}}

\newenvironment{proofof}[1]{\bigskip \noindent {\bf Proof of #1:}\ }{\qed}

\IEEEoverridecommandlockouts
\begin{document}

\title{On the Efficiency of Connection Charges---\\Part I: A Stochastic Framework}

\author{Daniel~Munoz-Alvarez
		and~Lang~Tong%<-this % stops a space
\thanks{This work is supported in part by the National Science Foundation under Grants CNS-1135844 and 15499.}% <-this % stops a space
\thanks{D. Munoz-Alvarez and L. Tong are with the School of Electrical and Computer Engineering, Cornell University, Ithaca, NY, 14853, USA.  Emails:  {\tt\small dm634@cornell.edu}, {\tt\small lt35@cornell.edu}}
}

\markboth{}%
{}

\maketitle %\thispagestyle{empty} \pagestyle{empty}

\begin{abstract}%
%% Abstract Part I
This two-part paper addresses the design of retail electricity tariffs for distribution systems with distributed energy resources such as solar power and storage.
In particular, the optimal design of dynamic two-part tariffs for a regulated monopolistic retailer is considered, where the retailer faces exogenous wholesale electricity prices and fixed costs on the one hand and stochastic demands with inter-temporal price dependencies on the other.
Part I presents a general framework and analysis for revenue adequate retail tariffs with advanced notification, dynamic prices and uniform connection charges.
It is shown that
%when the retailer is regulated to break even, 
%when the retailer revenue adequacy is required,
the optimal two-part tariff consists of a dynamic price that may not match the expected wholesale price and a connection charge that distributes uniformly among all customers the retailer's fixed costs and a price-volume risk premium.
A sufficient condition for the optimality of the derived two-part tariff among the class of arbitrary ex-ante tariffs is obtained.
Numerical simulations quantify the substantial welfare gains that the optimal two-part tariff may bring compared to the optimal linear tariff (without connection charge).
Part II focuses on the impact of two-part tariffs on the integration of distributed energy resources.
%Previous:
%The retail electricity tariff design problem with distributed renewable and storage resources is analyzed.
%In particular, the optimal design of an ex-ante linear dynamic tariff from the perspective of a regulator is studied as well as the allowance of a uniform connection charge.
%An exogenous wholesale electricity market, a monopolistic retailer with a revenue sufficiency constraint, and a model of heterogeneous residential customers with inter-temporal stochastic demands are considered. 
%This paper is divided in two parts.
%Part I develops the framework and studies the problem without distributed renewable and storage resources.
%It is shown that, when the retailer is guaranteed to break even, the optimal ex-ante two-part tariff consists of an hourly price that matches the expected real-time price and a connection charge that distributes uniformly the revenue target and a price-volume risk premium among all customers.
%A sufficient condition for the optimality of the derived two-part tariff among the class of arbitrary tariffs is obtained.
%Numerical simulations suggest that while the optimal linear tariff can bring surplus loses of 4.8\% the revenue collected by a utility in NYC relative to the utility's (suboptimal) two-part flat tariff, the optimal two-part tariff may bring surplus gains of 8.1\%.
%Part II analyses the integration of said resources by the retailer in the distribution network and by customers behind the meter.
\end{abstract}

% No keywords for conference version
\begin{IEEEkeywords}
Retail tariff design, connection charges, dynamic pricing, distributed energy resources, optimal demand response.\\[-0.5cm]
\end{IEEEkeywords}

\section{Introduction} \label{sec:intro}

\ifdefined\IEEEPARstart
	\IEEEPARstart{T}{he}
\else
	The
\fi
electric power industry is experiencing an important transformation driven by disruptive innovation in distributed renewable generation and energy storage systems \cite{GTM15}.
A concern of this transformation is the impact of the inclining adoption of said distributed energy resources (DERs) on the financial viability of regulated distribution grid operators \cite{RMI15}.
In particular, under the restriction to volumetric and net-metering tariffs, the gradual decline in energy sales could compromise the ability of grid operators to recover their predominantly fixed operational and capital expenditures.
This could result in the need to increase retail prices further above wholesale electricity prices, thereby amplifying the entailed economic inefficiencies, inter-customer cross-subsidies, and incentives for DER adoption in a vicious circle.

This two-part paper aims to shed lights on the effectiveness of connection charges as a means to mitigate the negative impacts of the sustained adoption of DERs.
To this end, Part I develops a framework to analyze the efficiency of retail electricity tariffs set for a regulated retailer who serves a heterogeneous population of residential customers under demand and wholesale price uncertainties.

In particular, we are interested in two practical and fairly general ex-ante retail pricing models: a volumetric linear tariff and a two-part tariff consisting of a volumetric linear charge and a connection charge.
Our goal in Part I is to gain insights into the structure of the optimal revenue adequate linear and two-part tariffs that allow us to analyze, in Part II, the effects of integrating customer and retailer-owned DERs.

%\red{
%In this context, this series of two papers aims to analyze from basic principles the design of retail electricity tariffs with connection charges considering $(i)$ a retailer revenue sufficiency constraint, and $(ii)$ the integration of renewable energy resources and energy storage systems by both the retailer and its customers.
%To that end, we revisit in this paper the retail tariff design problem from the perspective of a regulated retailer and present the integration analysis in the second paper.
%We aim to shed lights on the effectiveness of a connection charge as a means to mitigate the negative impacts of behind-the-meter distributed energy resources.
%}

%\blue{In particular, we shed lights on the role of fixed connection charges in the design of retail tariffs that maximize the collective consumer welfare and on their impact on inequity concerns.
%Our focus is on tariffs that are fixed and announced with some time lag to customers and that are valid for a billing period with multiple consumption periods.
%A proper analysis of such ex-ante tariffs demands the use of a stochastic demand model, which then becomes the natural framework to model the integration of renewables and storage.}

%In our work, we analyze how fixed connection charges can mitigate the inefficiency of volumetric charges in the presence of renewables and storage integrations in the distribution grid.
%In the long-run, such tariffs should be updated to induce a socially optimal level of investment in storage systems and renewable distributed generation.

\subsection{Related Work}
%\red{(Part 2)
%The integration of distributed energy resources (DER) in distribution networks is a complex dynamic process that involves modeling the way millions of customers decide to adopt and use such technologies.
%%Perhaps as a means to rationalize such a collective decision making process, this phenomenon has been studied at different time scales.
%To the best of our knowledge, the first study of this process and the \textit{death spiral hypothesis} was conducted by Cai \etal \cite{CaiEtal13}.
%Therein, authors model and analyze empirically the long-run closed loop dynamics of residential PV solar integration under California's current increasing block rate (IBR) structure.
%Their model predicts that the upward pressure on the tariffs induced by the IBR structure is ``unlikely to have a significant impact on future PV uptake rates in the next 10 years.''
%They also explore the delaying effect that connection charges may have on PV adoption rates.}

There is a vast literature on efficient retail pricing of electricity, the economic foundations of which reside in the classical theory of peak-load pricing \cite{BrownSibley86}---known more recently as dynamic pricing \cite{JoskowWolfram12}. %, which date back at least to \cite{Bye29}.
%These theories, which date back at least to \cite{Bye29, Boiteux49}, have been compiled in numerous works (\eg see \cite{BrownSibley86, Wilson93, CrewEtAl95}).
%There is an extensive literature on retail electricity pricing.
%Its economic foundations reside in the classical theories of public utility and peak-load pricing, which, dating back at least to \cite{Bye29, Boiteux49}, have been compiled in numerous works (\eg see \cite{BrownSibley86, Wilson93, CrewEtAl95}).
Recent reviews of the subject along with a brief history of its adoption by regulators and electric utilities can be found in \cite{BraithwaitEtal07,LBNL16b}.
In this context, connection charges have been considered as a means to raise additional revenue to recover the predominantly fixed costs of electric utilities \cite{BrownSibley86}.
%In this matter, some authors \cite{BorensteinHolland05, Joskow&Tirole06} have studied real-time linear and two-part tariffs in competitive and regulated market settings where investments in costly real-time meters are possible.
%These works relate to our setting since renewable and storage technologies, not unlike smart meters, bring uncertain short-run value and entail certain upfront capital costs to either the retailer, consumers, or both.
%While the availability of such technologies affect the tariff design problem in the short-run, the investment decision must be studied in the long-run.
In the U.S., mild connection charges are prevalent with exceptions such as California, where the large investor-owned utilities have default volumetric residential tariffs with virtually no connection charges%
\footnote{PG\&E and SDG\&E have no connection charge whereas SCE's charges \$0.99/month.
While these utilities have a minimum bill of \$10/month or less, it is binding on extremely few customers, and thus practically irrelevant \cite{Borenstein14}.}
%This is due in part to reasons dating back to the California electricity crisis and state legislators' interest to keep low-income households' bills stable \cite{Borenstein14a}.
\cite{Borenstein14}.
%although the largest investor-owned utilities (IOU) in California currently have small or no connection charges, significant connection charges are prevalent everywhere else in the US.

In the last two decades, while the adoption of time-varying prices has been particularly slow in the U.S. \cite{JoskowWolfram12}, the advent of cheaper smart meters, small-scale renewable energy installations, battery storage technologies and home energy management tools has stimulated research in \emph{dynamic pricing} (see \cite{DengEtAl15, VardakasEtAl15} and references therein) as sophisticated technologies 
%data analytic tools can help 
can enable customers to react to price signals \cite{FaruquiEtal12}. %SaadEtAl12
%The latter is a notion that broadly refers to retail prices that are allowed to vary in time in the short-term.
Of particular interest is real-time pricing (RTP), a form of dynamic pricing widely known to be a critical feature of efficient electricity markets \cite{BorensteinEtal02}.
An overview of dynamic pricing and a recent analysis of its limited adoption in the U.S. are available in \cite{BorensteinEtal02} and \cite{JoskowWolfram12}, respectively.

Economic approaches to dynamic pricing often rely on functional demand models to characterize competitive equilibrium prices when smart meters become available to customers \cite{Borenstein05, BorensteinHolland05, JoskowTirole06}.%JoskowTirole07
The most relevant analysis is \cite{JoskowTirole06} where the socially optimal linear and two-part retail tariffs subject to a retailer revenue sufficiency constraint are derived.
Unlike our work, however, this analysis does not accommodate inter-temporal demand dependencies nor the integration of DERs.
%Therein, in a case where customers have smart meters yet face a flat tariff, their setting accommodates either demand and wholesale price uncertainties \emph{or} multiple independent consumption periods, but not both cases.
%\red{In contrast, our setting fully accommodates both features and considers ex-ante time-differentiated tariffs.}

%Most engineering approaches, on the other hand, focus on analyzing the interactions between load-serving entities, aggregators, and end-customers \cite{DengEtAl15, VardakasEtAl15}. %SaadEtAl12
Most engineering approaches, on the other hand, focus on analyzing demand response models in smart grids \cite{DengEtAl15, VardakasEtAl15}. %SaadEtAl12
These approaches often involve modeling customer behavior \cite{MohsenianLeon10, MohsenianWong10}, sometimes down to the appliance level \cite{DaryanianEtAl89, ShaoEtAl11, MathieuEtAl15, LiEtAl16}, by characterizing customers' response to certain pricing scheme.
Appliances modeled include thermostatically controlled loads (TCL) \cite{MathieuEtAl15, LiEtAl16}, electric vehicles \cite{ShaoEtAl11}, and batteries \cite{DaryanianEtAl89}. %XuTong16
%In this community, dynamic pricing is classified as an indirect or price-based DR program along with time-of-use, critical peak pricing, among other pricing schemes.
Other works focus on designing said pricing schemes to induce a desired behavior anticipating customers' response \cite{ChenEtal12, LiEtAl13, YangEtAl13, MengZeng13, BinaAhmadi15, JiaTong16a}.
%\cite{LiEtAl13, BinaAhmadi15}, often posing the problem in a game-theoretic setting \cite{ChenEtal12, YangEtAl13, MengZeng13, JiaTong16a}. %SamadiEtAl10, LiEtAl11, HuangEtAl12, DoostizadehGhasemi12, Tarasak11, GkatzikisEtAl13, TangEtAl14
%MuratoriRizzoni16 (TOU, CPP)
%Most approaches characterize customers' optimal response to certain retail tariff, and some approaches further design such tariff to optimize some private or social objective, \eg to maximize social welfare.
% Liyan's and other related work
%Numerous works have focused on retail electricity pricing from different perspectives including that of consumers, suppliers, retailers, and the social planner.
%The long-term approach in \cite{CaiEtal13} is in contrast with the more numerous works that have focused on understanding the problem as posed to the stakeholders in the short-term.
%This is the spirit of \cite{ChenEtAl12}, \cite{JiaTong16a}, and our work.
For example, the work in \cite{ChenEtal12} considers utility-maximizing consumers and a social-welfare-maximizing supplier that procures electricity in two steps (day-ahead and in real-time).
In a multiperiod and deterministic setting, the authors derive the socially optimal retail prices: a time-differentiated linear tariff.
%No revenue sufficiency constraint is considered either.
%They formulate the social planner's problem as a stochastic dynamic program, and characterize the real-time prices that induce the optimal solution to such problem in the absence of uncertainty.
%Although authors formulate the problem as a constrained stochastic dynamic program where the consumer model is based on multiple appliances (including thermostatically controlled loads and energy storage), they only provide a solution algorithm for the deterministic case (\ie without renewable energy, which is the only source of uncertainty they consider).
%Their consumer model, based on multiple appliances (including thermostatically controlled loads, or TCLs, and energy storage), is versatile yet deterministic.
%They provide results and a solution algorithm only for the deterministic case (\ie without renewable energy).
In a similar setting where only the real-time market and a stylized TCL model are considered, the work in \cite{JiaTong16a} derives optimal day-ahead retail prices while accommodating cost and demand uncertainty and an implicit retailer revenue requirement\footnote{The revenue requirement is incorporated indirectly into the formulation through a weighted social welfare objective.}.

\subsection{Summary of Results and Contributions}

The main contribution of Part I is the explicit characterization of the optimal revenue adequate ex-ante two-part tariff for a stochastic demand with inter-temporal dependencies.
The results in Part I lay the foundation for analyzing the welfare impacts of integrating DERs such as solar power and energy storage under different retail tariffs, which we address in Part II.
Here we apply the classical Ramsey pricing theory with extensions to accommodate the uncertainty and inter-temporal dependencies of demand that arise with the integration of behind-the-meter renewables and storage.
While economic literature has disregarded said extensions, engineering approaches to dynamic pricing and demand response have ignored the revenue adequacy objective of retail tariff design, which our work addresses explicitly.
%thereby ignoring design criteria such as the critical issue of revenue adequacy.
%have focused on the pricing problem but consistently ignored the more general tariff design problem, thereby ignoring design constraints such as the critical issue of revenue adequacy.
%This is an increasingly relevant gap in the retail pricing problem that our work addresses.
%On the other hand, engineering approaches to demand response have focused on the pricing problem but consistently ignored the more general tariff design problem, thereby ignoring design constraints such as the critical issue of revenue adequacy.
%In this context, we aim at closing the wide gap existing in the retail pricing problem between both academic communities.
In this context, there are no existing comparable studies in the open literature with the exception of a preliminary work in \cite{MunozTong16}.
To a large extent, our results are an examination of emerging issues in smart distribution systems through the lens of  classical economic results on the fundamental efficiency of two-part tariffs in stylized economic models \cite{BrownSibley86}.

The main results of this paper are as follows.
We consider the design of ex-ante retail tariffs from the perspective of a regulated retailer subject to a revenue sufficiency constraint.
The ex-ante tariffs considered here include traditional tariffs with long lag times as well as some of the more sophisticated tariffs that are being considered for smart distribution systems.
Examples include time-of-use tariffs, critical peak pricing, variable peak pricing, and real-time pricing \cite{FaruquiEtal12}.
The retailer considered in this paper is a regulated monopoly which, on one hand, serves heterogeneous residential customers with elastic demands.
%The retailer, on one hand, serves residential customers with heterogeneous elastic demands.
The demand model considered here is stochastic and captures inter-temporal price dependencies.
On the other hand, the retailer interfaces with an exogenous wholesale market with stochastic real-time prices.
We describe the models in our formulation in Section \ref{sec:model}.
%The models integrating this general setting are described in Section \ref{sec:model}.

% summarize the main results
Within this general setting, we characterize the structure of optimal linear and two-part tariffs in the presence of demand and wholesale price uncertainty in Section \ref{sec:retailTariffDesign}.
In particular, we show that the optimal ex-ante two-part tariff consists of a time-varying retail price that not always matches the expected wholesale price and a connection charge that allocates uniformly among all customers the retailer's fixed costs and risk-related costs caused by the ex-ante determination of the tariff.
We further show that the optimal volumetric tariff, referred hereafter as linear tariff, is characterized by a time-varying price markup---relative to the optimal two-part tariff's price---that depends on the retailer's fixed costs and the price elasticity of demand.
%In particular, we show that the optimal volumetric tariff, referred hereafter as linear tariff, is characterized by a time-varying markup relative to the expected wholesale price that is function of the \emph{retailer surplus} (RS) requirement and the price elasticity of demand.
%For the two-part tariff, which has a connection charge independent of each customer's consumption, we show that the optimal connection charge is the sum of the share of the RS requirement and a price-volume risk premium customers pay to the retailer for facing uncertainty in demand and wholesale prices.

We further compare the efficiency of the linear and two-part tariff in Section \ref{sec:optimality}.
Specifically, we present a parametric characterization of the \emph{social welfare} (SW) or total surplus and the \emph{consumer surplus} (CS) as a function of the retailer's fixed costs.
We show that the two-part tariff achieves the same SW regardless of the retailer's fixed costs.
For the linear tariff, in contrast, the SW decreases as the fixed costs increase, thus characterizing a trade-off between fixed costs and efficiency.
We also provide a sufficiency condition under which the two-part tariff is optimal among all ex-ante nonlinear tariffs.

We demonstrate the performance of the derived tariffs numerically using publicly available data from NYISO and the largest utility company in New York City in Section \ref{sec:example}.
Contingent on the deployment of enabling technologies and smart meters, our results estimate that the optimal day-ahead linear tariff could bring loses ($4.8\%$ of the utility's revenue) relative to the utility's suboptimal two-part flat tariff due to the lack of a connection charge.
The optimal day-ahead two-part tariff, on the other hand, could bring significant gains ($8.1\%$ of the utility's revenue).
From a societal perspective, these loses and gains manifest themselves as reductions and increments in electricity consumption, respectively.
%Our results estimate that the optimal day-ahead two-part tariff can bring daily surplus gains at least $5$ times larger than the optimal day-ahead linear tariff.
%These gains, which represent dollar amounts within $6.5\%$ and $20\%$ of the revenue collected by the utility's flat rate, are likely to increase with the adoption of behind-the-meter DERs.
These estimates assume a realistic own-price elasticity of demand and a stylized model for TCL.
Some concluding remarks and proofs are included in Section \ref{sec:conclusions} and the Appendix, respectively.

\subsection{Notations}

%The paper is organized as follows.
%In Section \ref{sec:model}, we describe the retail tariff, consumer, and retailer models used in our formulation.
%This is followed by a specification of the tariff design problem and our main results in Section \ref{sec:retailTariffDesign}.
%Proofs of the main results are provided in the Appendix.
%In Section \ref{sec:illustrativeExample}, we offer a numerical example that demonstrates the efficiency gains brought the optimal linear and two-part tariffs.
%We close with some concluding remarks in Section \ref{sec:conclusions}.

% Exampel of similar way to incorporate proofs separately in TPS paper:
% From Baldick (2007) Border transmission rights
% Further details, additional discussion including comparison to other types of transmission rights, and technical results are established in a technical reference available online [20].

%As general notation, 
We use $\overline{x} = \Embb[x]$ to denote the expectation of a random vector $x \in \Rmbb^n$ and $\Sigma_{x,y} = \cov(x, y) \in \Rmbb^{n \times m}$ to denote the cross-covariance matrix of two random vectors $x \in\Rmbb^n$, $y \in \Rmbb^m$.
Let also $x_k$ denote the $k^{\text{th}}$ entry of a vector $x \in \Rmbb^n$ and $x^{\Top}$ its transpose.

\section{Model} \label{sec:model}

Given our focus on the retail electricity market, we assume the state of the wholesale market is represented by an exogenous discrete-time random process $\lambda_{k} \in \Rmbb_+$, which represents the wholesale RTP of electricity at time $k$ in a single location of interest.
We assume that the time periods $k=1,\ldots,N$ partition a billing cycle, which is the time horizon relevant for our formulation.
Moreover, we assume the wholesale RTP accurately reflects the \emph{social} marginal cost of electricity \cite{Borenstein14}.
%\footnote{As argued in \cite{Borenstein16}, the marginal cost of electricity may be lower than the social marginal cost due to the cost of negative externalities not beared by suppliers such as pollution. This requires us to assume, without loss of generality, that the RTP reflects the social marginal cost of electricity.}.

\subsection{A Retail Tariff Model}

In this paper, we consider time-differentiated retail electricity tariffs that are set and announced in advance (\ie ex-ante) by a regulated retailer with a fixed lag time.
These tariffs $(i)$ are fixed before the beginning of a billing period of certain length (\eg a month or a day) with a fixed lag time (\eg several days or hours), $(ii)$ specify a pricing rule that depends on the temporal consumption profile within the billing period rather than on the accumulated consumption, and $(iii)$ are allowed to vary dynamically from one billing period to the next.
In the context of retail tariffs, the tariff lag time induces a tradeoff between advanced price notification and price signal accuracy \cite{BorensteinEtal02}.
%While traditionally static electricity tariffs tend to have longer lag times (\eg years or months), tariffs based on dynamic pricing schemes usually exhibit shorter lag times (\eg days, hours, minutes, or even no advanced notification).
The tariff model considered here captures both the traditional long term flat tariff that has months or years of lag time as well as more sophisticated dynamic tariffs such as those with days or hours of advanced notification, but it generally excludes \textit{ex-post} tariffs such as those indexed to the wholesale RTP.
%While for end-customers shorter lag times often means more volatile prices, for retailers---who play an intermediary role between the wholesale and the retail market---they represent a way to share the price spread risks they face in a way that improves the overall social welfare.

%In this paper, we focus on a day-ahead hourly pricing tariff with and without a connection charge.
%The former is a type of dynamic pricing in which the retailer fixes an hourly electricity price and a connection charge one day ahead, but allows one or both parts to vary dynamically from one day to the next.
%In the context of dynamic pricing, such tariff provides an intermediate balance for customers between lower price volatility, greater advanced price notification, and more accurate price signals \cite{BorensteinEtal02}.
%Similarly, it allows retailers to partially share with their customers the price risks they face in the market.

%In practice, such a setting would allow the retailer to set the retail price based on its position in the forward and day ahead wholesale markets, and thus to share some price risk with its customers.
%More importantly, such tariff allows end consumers to respond to prices that reflect more accurately the short-run cost of producing electricity while still giving them some allowance to plan ahead \cite{BorensteinEtal02}.

Formally, some time before the billing cycle starts, the retailer announces a tariff $T:\Rmbb^N \rightarrow \Rmbb$ that maps the metered consumption power profile $q \in \Rmbb^N$ of each customer to a scalar charge $T(q) \in \Rmbb$.
While the $k^{\text{th}}$ entry of $q$ is a single customer's metered consumption in period $k$ of the billing period, the amount $T(q)$ (in dollars) represents the total bill.
%For each customer, the consumption $p$ represents his actual consumption.
%This implicitly assumes that customers do not have access to self-generation or storage technologies, which would allow them to offset and buffer their consumption of electricity from the power grid.
%Thus, in the absence of such behind-the-meter DERs the net-metered consumption clearly matches the actual consumption, \ie $p=q$.
Note that this form of tariff captures the intertemporal dependencies of pricing and consumption within each billing cycle (but not between several billing cycles).

Given a tariff $T$, customers rationally choose in real-time how much electricity to purchase from the retailer during each consumption period of the current billing cycle.
%In particular, they choose how much of that electricity is withdrawn from their own behind-the-meter renewables and storage rather than from the grid.
%We assume customers choose the amount they withdraw from the grid and their storage knowing the renewable generation available for their immediate consumption.
%In general, we reserve letter $p$ for \textit{net-metered} energy withdraws from the grid and $u$ for actual \textit{consumption} quantities.
The retailer then pays for the aggregate demand at the wholesale RTP.
%This amounts to paying for such net-metered consumption profile at the real-time wholesale price of electricity.

Although in practice retailers buy certain portions of the aggregate demand in forward markets (including the day-ahead market), we can neglect such purchases in our formulation without loss of generality for the following reason.
In perfectly competitive and well-functioning two-settlement markets, forward transactions are essentially used to hedge against the volatility of the RTP. %\cite[Chap. 3.2]{Stoft02}.
Here, we consider risk-neutral decision markers that deal with uncertainty by taking expectations.
Thus, in our setting, forward markets would bring no significant advantages to any stakeholder.
%In particular, they still face the the same performance incentives in the RT market as if they needed to trade all of their energy requirements at the real-time price \cite[Sec. 3-2.1]{Stoft02}.
This justifies the reliance of the retailer in the RTP to purchase electricity, which is an assumption that simplifies our exposition considerably.

\subsection{Consumer Model} \label{sec:consumer model}

We consider $M$ customers (indexed by $i$) who obtain a monetary benefit (\ie gross surplus) $S^i(q^i, \omega^i) \in \Rmbb$ from consuming a power profile $q^i \in \Rmbb^N$ throughout the billing cycle.
This benefit is contingent on $\omega^i=(\omega^i_1,\ldots,\omega^i_N) \in \Rmbb^N$, where $\{\omega^i_k\}_{k=1}^N$ is an exogenous random process that represents customer $i$'s local state.
%This benefit is state-contingent in the sense that $\omega^i=(\omega^i_1,\ldots,\omega^i_N) \in \Rmbb^N$ represents a random process whose $k^{\text{th}}$ entry denotes customer $i$'s exogenous state at time period $k$.
We assume that $S^i$ is continuously differentiable in $(q^i,\omega^i)$.
%We assume that $S^i$ is positive, increasing, strictly concave in $q^i$, and continuously differentiable in $q^i$ and $b^i$.

Accordingly, customer $i$ exhibits a consumption profile
$ q^i = q^i(T, \omega^i)$
when facing a tariff $T$ and a sequence of local states $\{\omega^i_k\}_{k=1}^N$.
Customers are rational in that sense that the sequence of consumptions $\{q^i_k(T,(\omega_1,\ldots,\omega_k))\}_{k=1}^N$ solves the multistage stochastic program
\begin{align} \label{eq:consumer surplus individual}
\overline{\cs}^i(T) = \max_{q^i(\cdot)} \ \Embb[ S^i(q^i(\omega^i), \omega^i) - T(q^i(\omega^i)) ],
\end{align}
where the expectation is taken over $\omega^i$, and $\overline{\cs}^i(T)$ represents customer $i$'s expected surplus.
%We emphasize here that $q^i(T,\omega^i)$, which is the solution of a multistage stochastic program, has a causal dependence on $\omega^i$.
%\red{We thus assume that the corresponding first-order necessary conditions
%\begin{align} \label{eq:FOC consumer problem}
%\Embb_{|k}[v^i_k(D^i(\pi,b^i),b^i)] = \pi_k, \quad k=1,\ldots,N 
%\end{align}
%%$$ \Embb_{|k}[v^i_k(q^i(T,b^i),b^i)] = \Embb_{|k}[\partial T(q^i(T,b^i))/ \partial q^i_k,$$
%are satisfied, where $v^i_k = \partial S^i / \partial q^i_k$ denotes the marginal benefit of consuming an incremental unit at time $t$, and the expectation $\Embb_{|k}$ is with respect to $b^i$ conditioned on $b^i_1,\ldots,b^i_k$.}
Correspondingly, a tariff $T$ yields an (aggregate) expected consumer surplus
\begin{align} \label{eq:cs}
\overline{\cs}(T) \ = \ \Embb\left[ \mbox{$\sum_{i=1}^M$} S^i(q^{i}(T,\omega^i), \omega^i) -  T(q^{i}(T,\omega^i)) \right],
\end{align}
where the expectation is taken over $\omega = (\omega^1,\ldots,\omega^M)$.
%where $q^i(T)$ denotes the metered consumption induced by tariff $T$ on customer $i$. 
%For notational convenience, we define the aggregate demand function $D(\pi) := \sum_{i=1}^M D^i(\pi,b^i)$\footnote{We drop the explicit dependence of $D(\pi)$ on each $b^i$ for notational convenience.}.
%For notational convenience, we define $D(\pi) := \sum_{i=1}^M D^i(\pi,b^i)$ and $S(q^1,\ldots,q^M) := \sum_{i=1}^M S^i(q^i,b^i)$ to refer to the aggregate demand and gross surplus\footnote{We drop the explicit dependence of $D(\pi)$ and $S$ on $b^i$ and $p^i$ only for notational convenience.}.
%A consumer model similar to the one decribed in this section for $N=1$ is proposed and discussed with more detail in \cite{JoskowTirole06}.
%More details of this demand response model can be found in the appendix and also in \cite{JiaTong16b, JiaTong13}.

Of particular interest is the demand response to tariffs $T$ with constant gradient $ \nabla T = \pi \in \Rmbb^N$, where $\pi \in \Rmbb^{N}$ is a time-varying per-unit price, such as the tariff with the affine form $T(q^i)=A+\pi^{\Top}q^i$.
For such tariffs $T$ we use the notation
\begin{align}
D^i(\pi, \omega^i) = q^{i}(T,\omega^i) \label{eq:demand function}
\end{align}
for customer $i$'s demand profile, thus implicitly assuming that it depends on $T$ only through $\pi$.
Hence, $D^i$ is a standard demand function which we assume to be nonnegative
%\footnote{\red{We relax the nonnegativity assumption on metered demand in the second part of this work to incorporate behind-the-meter renewables and storage.}}
and continuously differentiable in $\pi$, and its Jacobian $\nabla_{\pi} D^i \in \Rmbb^{N \times N}$, with $(k,t)$ entry $\partial D^i_k/\partial \pi_t$, negative definite.
%\red{For simplicity, we further assume that $\nabla D^i(\pi,b^i)$ and $\lambda$ are uncorrelated.}
Under the regularity assumptions made on $S^i$ and $D^i$, one can show that $\overline{\cs}^i(T)$ is decreasing and convex in $\pi$ (see Prop. \ref{prop:convexity cs} in Appendix).
%In Fig. \ref{fig:linear tariff}, we depict a demand function $D^i$ for the single period case ($N=1$) in the \textit{quantity-price} plane.
%Similarly, in Fig. \ref{fig:nonlinear tariff}, we illustrate the relationship between the demand function, the derivative $T'$ of a nonlinear tariff, and the induced consumer surpluses in the quantity-price plane for the single period case ($N=1$).
We further define the aggregate demand function as $D(\pi,\omega) := \sum_{i=1}^M D^i(\pi,\omega^i)$.
A consumer model similar to the one decribed in this section for $N=1$ is proposed and discussed with more detail in \cite{JoskowTirole06}.

For example, for a linear tariff $T(q^i) = \pi^{\Top} q^i $, the consumption of a TCL may be modeled with a linear demand function $D^i(\pi,\omega^i)=\omega^i-G^i\pi$, with deterministic and positive definite $G^i \in \Rmbb^{N \times N}$.
Such demand function can be derived from an additive and temporally-separable quadratic benefit function $S^i$ via stochastic dynamic programming \cite{JiaTong16a}. %JiaTong12 (Allerton),JiaTong13 (CDC), JiaTong15 (PESGM) left out

\subsection{Retailer Model} \label{sec:retailer model}

We consider the case of a retail monopoly and refer to the single entity as the retailer, utility, or load-serving entity (LSE).
In procuring an aggregate demand profile $q = \sum_{i=1}^M q^i \in \Rmbb^N$, we assume that the retailer incurs a variable cost $\lambda^{\Top} q$, where $\lambda = (\lambda_1,\ldots,\lambda_N)^{\Top} \in \Rmbb^N$ is the wholesale RTP\footnote{While in practice the time granularity of the wholesale RTP is finer than that of retail rates, we assume they are equal to simplify our exposition.
%While in practice the time scale of the RTP is typically on the minutes level (\eg every 5 or 15 minutes), we refer here with RTP to the average price associated to each interval $k=1,\ldots,N$ of the billing period. For example, if day-ahead hourly pricing is considered, $\lambda_k$ would correspond to the RT price averaged across the $k$-th hour of the day.
}
%the actual settlements often use average hourly prices presumably as a way to match the time scale of day-ahead prices.}.
%For each time $k$, we model the latter as a state-dependent nonnegative random variable $\lambda_k(\xi_k)$.
Hence, a tariff $T$ yields the expected retailer surplus
\begin{align} \label{eq:rp}
\overline{\rp}(T) \ = \ \Embb \left[ \mbox{$\sum_{i=1}^{M}$} T(q^{i}(T,\omega^i)) \ - \lambda^{\Top} q^i(T,\omega^i) \right],
\end{align}
where the expectation is taken over the global state $\xi = (\lambda,\omega)$.
For notational convenience, we define the RS collected from the volumetric charge $\pi$ of an affine tariff $T(q)=A+\pi^{\Top}q$ as $\phi(\pi) = (\pi - \lambda)^{\Top}D(\pi,\omega)$ so that $\overline{\rp}(T) = \overline{\phi}(\pi) + M A$ and
\begin{align} \label{eq:phi}
\overline{\phi}(\pi) 	&= (\pi - \overline{\lambda})^{\Top}\Embb[D(\pi,\omega)] - \tr(\cov(\lambda,D(\pi,\omega))).
\end{align}

\section{Retail Tariff Design} \label{sec:retailTariffDesign}

In our retail tariff design framework, we assume that the regulator mandates the retailer to choose a tariff $T$ that maximizes the expected consumer surplus.
Moreover, in order to recover the upstream fixed costs incurred to deliver electricity, the tariff should satisfy the revenue adequacy constraint $\overline{\rp}(T) = F$, where $F$ is a target approved by the regulator.

Formally, the regulator's problem can be stated as
\begin{align} \label{eq:reg problem}
\max_{T(\cdot)}	&	\	\overline{\cs}(T) \quad \text{s.t.} \quad \overline{\rp}(T) = F.
\end{align}
In general, this problem falls in the category of \emph{Ramsey-Boiteux pricing} and \emph{peak-load pricing} in economics \ifThesis{\cite{Ramsey27, Boiteux71, CrewEtAl95}}{\cite{CrewEtAl95}}, which are main components of the theory of public utility pricing \cite[Sec. 4.5]{BrownSibley86}.
See \cite{LBNL16b} for a recent overview of this problem in the context of electricity pricing.

In this section, we study linear and two-part tariffs---two of the most widely used tariffs in electricity industry.
The restriction to two-part tariff can in fact be made without loss of generality under certain conditions (Theorem \ref{thm:optimality}).
We begin by making the following assumption that guarantees the existence and uniqueness of solutions to problem \eqref{eq:reg problem}.

\assumptionAlt{ \label{Assumption 1}
$g(\pi)=\Embb[\nabla_{\pi} D(\pi,\omega) (\pi-\lambda)]$ is such that the Jacobian matrix $\nabla g(\pi)$ is negative definite (nd).
%is \emph{decreasing} in $\pi$ in the sense that $(g(\pi) - g(\tilde{\pi}))^{\Top} (\pi - \tilde{\pi}) \leq 0 $, $\forall$ $\pi$ and $ \tilde{\pi}$.
}
%\noindent This technical assumption, also made in \cite{JoskowTirole06} for the single period case, has limited practical implications\footnote{For instance, in \cite{JoskowTirole06}, authors note that this assumption is satisfied when $N=1$ if the curvature of the demand function is small enough (\ie $|D''(\cdot,\omega)/D'(\cdot,\omega)|$ small).}.
%\red{
%Intuitively, this assumption and the negative definiteness of $\Embb[\nabla_{\pi} D(\pi,\omega)]$ are necessary and sufficient for the (strict) concavity (in $\pi$) of the expected retailer surplus, $\overline{\rp}(T)$, under an affine tariff $T$.
%}
%%Note it is satisfied by the linear demand function since $G=\sum_{i=1}^M G^i$ must be negative definite.}.
%%We use this typical assumption, also made in \cite{JoskowTirole06}, to guarantee the existence of a solution to \eqref{eq:reg problem}.
%
%\assumptionAlt{ \label{Assumption 1}
%The expected retailer surplus $\overline{\rp}(T)$ is strictly concave in the prices $\pi$ of any affine tariff $T(q)=A+\pi^{\Top}q$.
%}
\noindent This assumption---made mainly for analytical convenience---is common is economics \cite{JoskowTirole06} and essentially imposes a limitation on the curvature of the demand function\footnote{%
%%***********************************************************************
%% Replaced by footnote below due to change in statement of Assumption.
%%***********************************************************************
%Notably, A\ref{Assumption 1} is implied by the negative definiteness (nd) of the Jacobian $\nabla_{\pi} D^i$ for a linear demand since the Hessian $\nabla^2_{\pi} \overline{\rp}(T) = 2\Embb[\nabla_{\pi} D]$ is nd.
%Moreover, for the case of additive disturbances, \ie $D(\pi,\omega)=\omega+D(\pi)$, the concavity of each function $D_{k}(\pi)$ in $\pi$ suffices because $\nabla^2_{\pi} \overline{\rp}(T) = 2\nabla_{\pi} D + \sum_{k=1}^N (\pi_k-\overline{\lambda}_k)\nabla_{\pi}^2 D_k$ is concave for all $\pi \geq \overline{\lambda}$.
%When any such demand is strictly convex in $\pi$, however, $\nabla^2_{\pi} \overline{\rp}(T)$ must be nd for A\ref{Assumption 1} to hold.
%It thus suffices the curvature of each demand, $|\nabla_{\pi} D^{-1} \nabla^2_{\pi} D_k|$, to be small.
%In \cite{JoskowTirole06}, for example, authors suggest that A\ref{Assumption 1} is satisfied when $N=1$ if the curvature $|D''(\cdot,\omega)/D'(\cdot,\omega)|$ is small enough.
%%*************************************************************************
In particular, for a linear demand $D(\pi,\omega)=b(\omega)-G(\omega)\pi$, $\nabla g(\pi)=\Embb[\nabla_{\pi} D(\pi,\omega)]$ is nd since $\nabla_{\pi} D(\pi,\omega)$ is nd.
Moreover, for a demand with additive disturbances $D(\pi,\omega)=b(\omega) + D(\pi)$, A\ref{Assumption 1} holds for $\pi \geq \overline{\lambda}$ if each $D_{k}(\pi)$ is concave in $\pi$ since $\nabla g(\pi) = \nabla D(\pi) + \sum_{k=1}^N (\pi_k-\overline{\lambda}_k) \nabla^2 D_k(\pi)$.
Concave demand functions are common in economic models since they guarantee profit and welfare maximization problems to be well defined \cite{AguirreEtal10}. %Malueg94
See, for example, Prop. \ref{prop:concavity rs} in the Appendix.
%In particular, A\ref{Assumption 1} implies the strict concavity of $\overline{\rp}(T)$ and of the social welfare metric $\overline{\sw}(T) = \overline{\cs}(T)+\overline{\rp}(T)$ in the prices $\pi$ of affine tariffs $T(q)=A+\pi^{\Top}q$ (see Prop. \ref{property:concavity rs} in the Appendix).
%, thus facilitating their optimization.
}.
Intuitively, the demand can be generally linear, concave, or convex in $\pi$; however, when convex, restrictions on the ``amount'' of convexity are required for Assumption \ref{Assumption 1} to hold.

\subsection{Structure of Optimal Two-Part Tariff} \label{sec:two part tariff}

By restricting the regulator's problem to two-part tariffs of the form $T(q)=A+\pi^{\Top} q$, problem \eqref{eq:reg problem} can be reformulated as a convex program under Assumption \ref{Assumption 1}.
We emphasize here that our analysis implicitly assumes that no customer chooses to avoid the connection charge by not consuming electricity at all\footnote{This assumption is widely accepted for services such as electricity and water since ``it is extremely unlikely that a customer will drop out of the market, however high the tariff'' \cite[Sec. 4.5]{BrownSibley86}.
Studies suggest, however, that more cost-effective DERs might challenge this assumption in future years \cite{RMI15}.}.
The following result characterizes the optimal solution.

\theoremAlt[(Optimal two-part tariff)]{\label{thm:affineTariff}
The two-part tariff $T^{*}$ that solves problem \eqref{eq:reg problem} is characterized by
\begin{align} 
%\pi^* &= \red{\Embb[\nabla_{\pi}D(\pi^*,\omega)]^{-1} \Embb[\nabla_{\pi}D(\pi^*,\omega) \lambda],} \label{eq:affineTariff2} \\[5pt]
\pi^* &= \overline{\lambda} + \Embb[\nabla_{\pi}D(\pi^*,\omega)]^{-1} \Embb[\nabla_{\pi}D(\pi^*,\omega) (\lambda-\overline{\lambda})], \label{eq:affineTariff} \\[5pt]
A^* &= \mbox{$\frac{1}{M}$} \left( F - \overline{\phi}(\pi^*) \right). \label{eq:connection charge}
%\phi(\pi) &:=  \Embb[(\pi^* - \lambda)^{\Top}D(\pi^*,\omega)]
%\Embb[\nabla_{\pi}D(\pi^*,\omega) (\pi^*-\lambda) ] = 0
\end{align}
}
Theorem \ref{thm:affineTariff} implies that the optimal price $\pi^*$ is characterized by a \emph{period-specific} price markup relative to the expected RTP, $\overline{\lambda}$.
Examination of \eqref{eq:affineTariff} reveals that this markup is essentially determined by the cross-covariance between the price sensitivity of demand and the RTP\footnote{In expression \eqref{eq:affineTariff}, the second expectation is a second-order expectation that can be thought as the cross-covariance between a matrix and a vector.}
To gain intuition into \eqref{eq:affineTariff}, consider a demand independent across time\footnote{That is, a demand with $D_k(\pi,\omega)$ independent of $\pi_t$ for all $t\neq k$}, case in which
\begin{align} \label{eq:affineTariff2}
\pi^*_k = \overline{\lambda}_k + \overline{\varepsilon}_{kk}(\pi^*)^{-1} \cov \left(\varepsilon_{kk}(\pi^*),\lambda_k \right),
\end{align}
for each $k=1,\ldots,N$, where we use
\begin{align} \label{eq:elasticity}
\varepsilon_{kt}(\pi) =\frac{\partial D_k(\pi,\omega) / \partial \pi_t}{\Embb[D_k(\pi,\omega)]/\pi_t}
\end{align}
to represent the (own or cross-time) price elasticity of demand at time $k$ with respect to the price at time $t$.
The latter result resembles the (second-best optimal) ex-ante two-part tariff derived in \cite[Sec. 3]{JoskowTirole06} for the single period case ($N=1$)\footnote{Which also applies to the continuous time case where $k\in[0,1]$, the demand and prices are deterministic, and demand is independent across time.}.

The expression \eqref{eq:connection charge} for the optimal connection charge $A^*$ also has an intuitive interpretation.
The first term corresponds to a \emph{uniform contribution} towards the retailer's target $F$.	
And, the second term corresponds to a \emph{uniform preallocation} of the surplus that the retailer expects to collect from the volumetric charge $\pi^*$, $\overline{\phi}(\pi^*)$, which---as noticeable from \eqref{eq:phi}---may be positive or negative in general.

To gain additional insights into these results, we have the following corollary.

\corollaryAlt{\label{cor:affineTariff1}
If $\nabla_{\pi} D(\pi,\omega)$ and $\lambda$ are uncorrelated, then $\pi^* = \overline{\lambda}$ and $A^*=\mbox{$\frac{1}{M}$} \left( F + \tr(\cov(\lambda, D(\overline{\lambda},\omega))) \right)$.
}

Corollary \ref{cor:affineTariff1} indicates that $T^*$ has a very simple and appealing structure that resembles the result for the \emph{deterministic} case where $\pi^*=\lambda$ and $A^*=F/M$.
Note that the assumption made in Corollary \ref{cor:affineTariff1} is valid for many situations.
It is certainly true for demands that are not much affected by consumers' local randomness\footnote{More precisely, this assumption is satisfied by demands whose \emph{sensitivity} to prices depends on the customers' set of appliances and idiosyncratic preferences rather than on random exogenous factors affecting the wholesale prices, such as random temperature fluctuations.}, such as the charging of electric vehicles and typical household appliances.
Even for loads from smart HVAC systems that are affected by random temperature fluctuations, the assumption in Corollary \ref{cor:affineTariff1} holds because the demand function takes the form $D(\pi,\omega) = \omega + D(\pi)$ \cite{JiaTong16a}.

As for the simpler structure of $T^*$, it may not be surprising since the efficiency of marginal cost pricing (\ie $\pi^* = \overline{\lambda}$) is a classical result for the deterministic case \cite[Sec. 4.5]{BrownSibley86}\cite{JoskowTirole06}.
Intuitively, marginal cost pricing is efficient because it induces customers to increase consumption until the derived marginal benefit matches the marginal cost of procuring electricity.

The expression for $A^*$ in Corollary \ref{cor:affineTariff1} also has an intuitive interpretation.
While the first term remains unchanged from \eqref{eq:connection charge}, the second term becomes a risk premium associated to the cross-correlation that the demand and the RTP may exhibit.
When such cross-correlation is positive (as in practice \cite{Burgeretal07}), the retailer is likely to face additional variable costs since the expected variable cost $\Embb[\lambda^{\Top}D(\overline{\lambda},\omega)]$ is larger than the variable revenue $\overline{\lambda}^{\Top} \Embb[D(\overline{\lambda},\omega)]$.
Intuitively, this fee represents a \emph{uniform risk premium} that customers pay to face a deterministic price rather than the volatile RTP.
Presumably, the inter-customer cross-subsidies arising from the uniform allocation of this risk premium  are negligible compared to the differences.
However, the integration of behind-the-meter renewables could make these cross-subsidies worth adjusting, for example, through the use of discriminatory connection charges consistent with the cost-causation principle described in \cite{PerezSmeers03}.
A discussion on cross-subsidies is held in Part II \cite{MunozTong16partII}.

\subsection{Structure of the Optimal Linear Tariff} \label{sec: linear tariff}

A tariff of the form $T(q^i)=\pi^{\Top} q^i$---a linear tariff---is an ex-ante two-part tariff with no connection charge.
While such purely volumetric tariff may be simpler, it has two fundamental disadvantages.
First, a closed form expression of the optimal linear tariff is not available under general assumptions.
Second, such restriction introduces a fundamental trade-off between the retailer surplus target and the attainable social welfare.
These drawbacks are noticeable in Theorem \ref{thm:linearTariff} and Corollary \ref{cor:linearTariff}, respectively.

When restricted to linear tariffs, a unique solution to problem \eqref{eq:reg problem} can be obtained due to Assumption \ref{Assumption 1}.
We characterize the optimal solution in the following result.

\theoremAlt[(Optimal linear tariff)]{\label{thm:linearTariff}
Consider the regime where $F$ is large, \ie $F \geq \overline{\phi}(\pi^*)$.
If feasible, the linear tariff $T^{\dagger}$ that solves problem \eqref{eq:reg problem} is characterized by
\begin{align} \label{eq:linearTariff}
\pi^{\dagger} = \pi^{*} -  \mbox{$\frac{\gamma-1}{\gamma}$} \Embb[\nabla_{\pi} D(\pi^{\dagger},\omega)]^{-1} \Embb[ D(\pi^{\dagger},\omega)],
\end{align}
or, equivalently, by
\begin{align} \label{eq:linearTariffb}
\sum_{t=1}^{N} -\overline{\varepsilon}_{kt}(\pi^{\dagger}) \left( \frac{\pi^{\dagger}_t - \pi^*_t}{\pi^{\dagger}_t} \right)  = \mbox{$\frac{\gamma-1}{\gamma}$}, \ \ \ \forall \ k=1,\ldots,N, \\[-20pt] \nn
\end{align}
where $\gamma$, the Lagrange multiplier of \eqref{eq:reg problem}, satisfies $\mbox{$\frac{\gamma-1}{\gamma}$} \in [0,1]$ and is such that $\overline{\rp}(T^{\dagger}) = F $.
In this regime of $F$, the problem is feasible if and only if
$F \leq \overline{\phi}(\pi^{\textsc{m}})$,
where $\pi^{\textsc{m}}$ is the price that maximizes $\overline{\rp}(T)$ over $\pi$, which matches $\pi^{\dagger}$ as $\gamma \rightarrow \infty$.
}

In Theorem \ref{thm:linearTariff}, expression \eqref{eq:linearTariff} reveals that the structure of the optimal linear tariff is characterized by a \emph{period-specific} price markup relative to the price of the optimal two-part tariff $\pi^*$.
The scalar $\frac{\gamma-1}{\gamma}\in [0,1]$, often called the Ramsey number, adjusts markups in all periods \emph{uniformly} to the point where the expected retailer surplus matches the target $F$.
A closer examination of \eqref{eq:linearTariff}, which can be rewritten as \eqref{eq:linearTariffb}, shows that the own and cross price elasticities of demand determine altogether the markup for each period within the billing cycle.

To understand \eqref{eq:linearTariffb}, it is informative to consider the case where the demand is independent across time, namely, $\varepsilon_{kt}(\cdot)=0$ for $t\neq k$.
In this case, the product of the markup $(\pi^{\dagger}_k - \pi^*_k)/\pi^{\dagger}_k$ and the own-price elasticity $-\varepsilon_{kk}(\pi^*)$ remains constant in time and equal to the Ramsey number.
This means that periods with inelastic demands get high markups and periods with elastic demands get low markups. 
For this reason, this pricing rule is known in economics as the \emph{inverse elasticity rule} \cite[Sec. 3.3]{BrownSibley86}.

Even simpler is the single period case, also derived in \cite[Sec. 3]{JoskowTirole06}.
Notably, when $N=1$, the scalar price $\pi^{\dagger}$ can be obtained directly from the constraint $\overline{\rp}(T^{\dagger})=F$, and it must be set so that it pays for the average total cost of the procured electricity, \ie $\pi^{\dagger} = (\Embb[ \lambda \cdot D(\pi^{\dagger},\omega)] + F)/\Embb[D(\pi^{\dagger},\omega)]$.

A specialized application of this result was developed in \cite{JiaTong16a}, where a model for TCLs under a day-ahead hourly pricing scheme was considered.
In this case, the demand function for each consumer is linear and the surplus function is quadratic \cite{JiaTong16a}.
The aggregated demand is therefore also linear.
The consumers as a collective have a quadratic aggregated surplus \cite{JiaTong16a}. Specifically,
\begin{align}
D^i(\pi,\omega^i) &= \omega^i - G^i \pi, \label{eq:demand function linear} \\
S^i(D^i(\pi,\omega^i),\omega^i) &= \delta^i(\omega^i) - \mbox{$\frac{1}{2}$} \pi^{\Top} G^i \pi. \label{eq:benefit function quadratic}
\end{align}
where $G^i \in \Rmbb^{24 \times 24}$ is deterministic, positive definite (and symmetric).
Letting $G=\sum_{i=1}^M G^i$ and $\Omega=\sum_{i=1}^M \omega^i$ and applying Theorem \ref{thm:linearTariff} readily yields
\begin{align} \label{eq:linearTariff2}
\pi^{\dagger} = \overline{\lambda} + \mbox{$\frac{\rho}{1+\rho}$} \hspace{1pt} (\pi^{o} - \overline{\lambda}),
\end{align}
where $\pi^{o}=G^{-1}\overline{\Omega}$ induces $\Embb[D(\pi^{o},\Omega)]=0$ and $\rho = \mbox{$\frac{\gamma-1}{\gamma}$}$ is the Ramsey number, which is set so that $\overline{\rp}(T^{\dagger})=F$.
Intuitively, $\rho$ varies within $[0,1]$ inducing prices that vary between $\overline{\lambda}$ and the profit-maximizing price $\pi^{\textsc{m}}=\frac{1}{2}(\pi^{o}+\overline{\lambda})$ as $F$ varies between $\overline{\phi}(\overline{\lambda})$ and the maximum profit $\overline{\phi}(\pi^{\textsc{m}})$.

\subsection{Tariff Performance Comparison} \label{sec:optimality}

\setlength{\fboxrule}{0pt}%

\begin{figure}
\centering
\fbox{\includegraphics[width=0.7\linewidth]{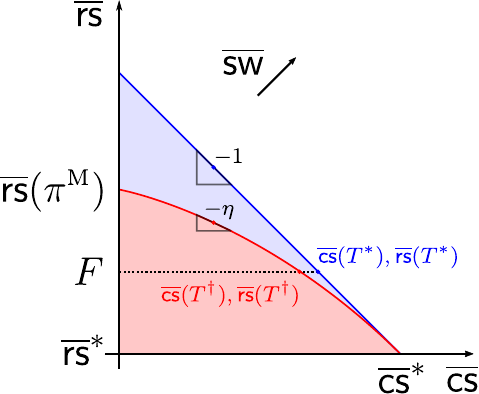}}
\caption{Pareto fronts induced by optimal linear and two-part tariffs.
The slope $\partial \overline{\rp}(F)/ \partial \overline{\cs}(F) = -\eta$ relates to $\gamma \geq 1$ in \eqref{eq:linearTariff} through $\eta=\frac{2}{1+\gamma}$.}
\label{fig:noIntegrationGeometry}
\end{figure}

We now discuss the performance of the derived tariffs in terms of social welfare (expected total surplus) leveraging the graphical representation provided in Fig. \ref{fig:noIntegrationGeometry}.
Therein, a Pareto front for each tariff illustrates the expected CS and SW induced by the tariff for different RS targets $F \in [\overline{\phi}(\pi^*),\overline{\phi}(\pi^{\textsc{m}})]$.
On one hand, Theorem \ref{thm:affineTariff} has the following implication.

\corollaryAlt{\label{cor:affineTariff2}
As a tariff parametrized by $F$, the two-part tariff $T^*$ induces an expected total surplus $\overline{\sw}^*$ that is independent of $F$ and $\overline{\cs}(T^*) =\overline{\sw}^* - F$.}

In Corollary \ref{cor:affineTariff2}, $\overline{\sw}^*$ denotes the constant SW attained by the price $\pi^*$, where $\overline{\sw}^*=\overline{\sw}(T^*)=\overline{\cs}(T^*)+\overline{\rp}(T^*)$.
Implicit in this result is that for any affine tariff one can check that
\begin{align} \label{eq:sw_optimal}
\overline{\sw}(T) = \mbox{$\sum_{i=1}^{M}$} \Embb \left[ S^i(D^{i}(\pi,\omega^i), \omega^i) - \lambda^{\Top} D^{i}(\pi,\omega^i) \right]
\end{align}
depends on $\pi$ but not on $A$.
Corollary \ref{cor:affineTariff2} thus implies that under the tariff $T^*$, collecting additional revenue from customers to cover larger fixed costs embedded in $F$ reduces consumer welfare but does not compromise social welfare.
This ``iso-efficient'' trade-off between retailer and consumer surplus is illustrated in Fig. \ref{fig:noIntegrationGeometry} with a \emph{linear} Pareto front with negative and unitary slope in the $\overline{\cs}$-$\overline{\rp}$ plane.
Intuitively, suboptimal two-part tariffs\footnote{As the ones currently used by most utilities in the U.S., due in part to additional bill stability and alleged equity concerns imposed by regulators.} can achieve any point in the $\overline{\cs}$-$\overline{\rp}$ plane in the shaded area below the linear Pareto front in Fig. \ref{fig:noIntegrationGeometry}, but no ex-ante two-part tariff can achieve points above this front.

Theorem \ref{thm:linearTariff}, on the other hand, has a analogous implication.
\corollaryAlt{\label{cor:linearTariff}
The quantities $\overline{\sw}(T^{\dagger})$ and $\overline{\cs}(T^{\dagger})$ induced by $T^{\dagger}$ as a tariff parametrized by $F$ are decreasing and concave in $F \in [\overline{\phi}(\pi^*),\overline{\phi}(\pi^{\textsc{m}})]$ with $\overline{\sw}(T^{\dagger}) = \overline{\sw}^*$ and $\overline{\cs}(T^{\dagger}) = \overline{\sw}^* - F$ for $F = \overline{\phi}(\pi^*)$.
}
Corollary \ref{cor:linearTariff} reveals that, unlike the tariff $T^*$, the optimal linear tariff $T^{\dagger}$ compromises not only consumer welfare but also social welfare when collecting additional revenue from customers is required to cover larger fixed costs embedded in $F$.
This trade-off is depicted in Fig. \ref{fig:noIntegrationGeometry} with a \emph{decreasing and concave} Pareto front in the $\overline{\cs}$-$\overline{\rp}$ plane that bends away from the efficiency level $\overline{\sw}^*$ attained by the tariff $T^*$ as $F$ increases from $\overline{\phi}(\pi^*)$ until it reaches $\overline{\phi}(\pi^{\textsc{m}})$.
As before, suboptimal linear tariffs can achieve any point in the shaded area below the curved Pareto front in Fig. \ref{fig:noIntegrationGeometry}, but no ex-ante linear tariff can achieve points above this front.

From the previous analysis, it is clear that two-part tariffs dominate linear tariffs in terms of expected consumer surplus in the regime of practical relevance where $F \geq \overline{\phi}(\pi^*)$.
A natural question to ask is whether two-part tariffs can be dominated by more complex nonlinear ex-ante tariffs.
We now argue that, under certain sufficient condition, the two-part tariff $T^{*}$ is  indeed optimal for the regulator's problem \eqref{eq:reg problem} among all ex-ante arbitrary tariffs.
To establish such result it suffices to show that $T^{*}$ induces the same expected consumer surplus that would be achieved by a social planner who makes consumption decisions on behalf of customers with the unconstrained objective of maximizing the expected total surplus.
This is because the social planner's problem provides a trivial upper bound to the regulator's problem.

Because we are interested in comparing ex-ante tariffs only, the social planner's problem should incorporate such implicit restriction.
The restriction to ex-ante tariffs translates into a restriction for the social planner to use only the information observable by each customer $i$ when choosing their consumption, namely their local state $\omega^i$.
Hence, the social planner's problem can be stated as
\begin{align} \label{eq:social planner}
\max_{\{q^i(\omega^i)\}_{i=1}^{M}} \overline{\sw} = \Embb \left[ \mbox{$\sum_{i=1}^{M}$} S^i(q^{i}(\omega^i), \omega^i) - \lambda^{\Top} q^{i}(\omega^i) \right],
\end{align}
where the expectation is taken with respect to $\xi=(\lambda,\omega)$, and $q^i(\omega^i)$ is causally contingent on (\ie adapted to) the local state $\omega^i$.
Finally, under the assumption that each $\omega^i$ and $\lambda$ are independent, we show that the optimal solution to \eqref{eq:social planner} is $q^i(\omega^i) = D^i(\pi^*,\omega^i)$, which matches the demand induced by the optimal two-part tariff.
This result and the implied optimality of $T^*$ are summarized in the following Theorem.

\theoremAlt{\label{thm:optimality}	% [(Optimality of two-part tariff)]
If \emph{(A2)} the wholesale RTP $\lambda$ and the local state $\omega^i$ of each customer $i=1,\ldots,M$ are statistically independent, then the two-part tariff $T^*$ is an optimal solution of \eqref{eq:reg problem} among all arbitrary tariffs with the same lag time.
}

Theorem \ref{thm:optimality} indicates that the restriction to two-part tariffs may imply no loss of generality.
This applies---less generally than Corollary \ref{cor:affineTariff1}---for demands that are not affected by consumers' local randomness\footnote{In other words, demands that, given a retail price vector $\pi$, are independent from random exogenous factors affecting the wholesale prices $\lambda$.}, such as washers and dryers, computers, and the charging of electric vehicles.
For all the other cases, where the sufficient condition (A2) does not hold, Theorem \ref{thm:optimality} sheds lights on the performance of the optimal ex-ante two-part tariff $T^*$.
While it is clear that the ex-ante restriction entails some efficiency loss when (A2) is not satisfied\footnote{This is because relaxing the ex-ante restriction enables the use of the ex-post two-part tariff $T(q^i)=F/M+ \lambda^{\Top}q^i$ which trivially achieves the maximum social welfare a social planner could achieve (first-best).}, it is not clear whether the restriction to two-part tariffs does entail efficiency losses.
In other words, is there a necessary condition for $T^*$ to be an optimal solution of \eqref{eq:reg problem}?

\section{Numerical Example} \label{sec:example}

\def\figwid{0.38\linewidth} % For 2 fig per row

\setlength{\fboxsep}{0pt}%
\setlength{\fboxrule}{0pt}%

%\begin{figure*}[t]
%\centering
%\hspace{10pt}
%\subfloat[][$F \in \left(\overline{\rp}(\pi^*),\overline{\rp}(\pi^\textsc{m}) \right)$.]{
%\label{fig:example:Pareto}
%\fbox{\includegraphics[width=\figwid, trim=41pt 0pt 0pt 35pt, clip=false]{figures/ParetoFront-Normalized-noDER-M2200000-alpha2e-001.eps}}
%} \hfill %\hspace{10pt}
%\subfloat[][Zoom into positive quadrant of \subref{fig:example:Pareto}.]{
%\label{fig:example:ParetoZoom}
%\fbox{\includegraphics[width=\figwid, trim=41pt 0pt 0pt 35pt, clip=false]{figures/ParetoFront-Normalized-noDER-M2200000-alpha2e-001-ZoomIn.eps}}
%}
%\hspace{10pt}
%\caption{CS and RS gains induced by various tariffs parametrized by the RS target $F$ assuming an intermediate own-price elasticity of $\overline{\varepsilon}(\pi^{\textsc{ce}})=-0.3$.}
%\label{fig:ParetoFronts}
%\end{figure*}

% Conditional to add figures for 1 or 2 column
\ifdefined\IEEEPARstart
% Two-column version here

\begin{figure*}[t]
\centering
\includegraphics[width=1\linewidth]{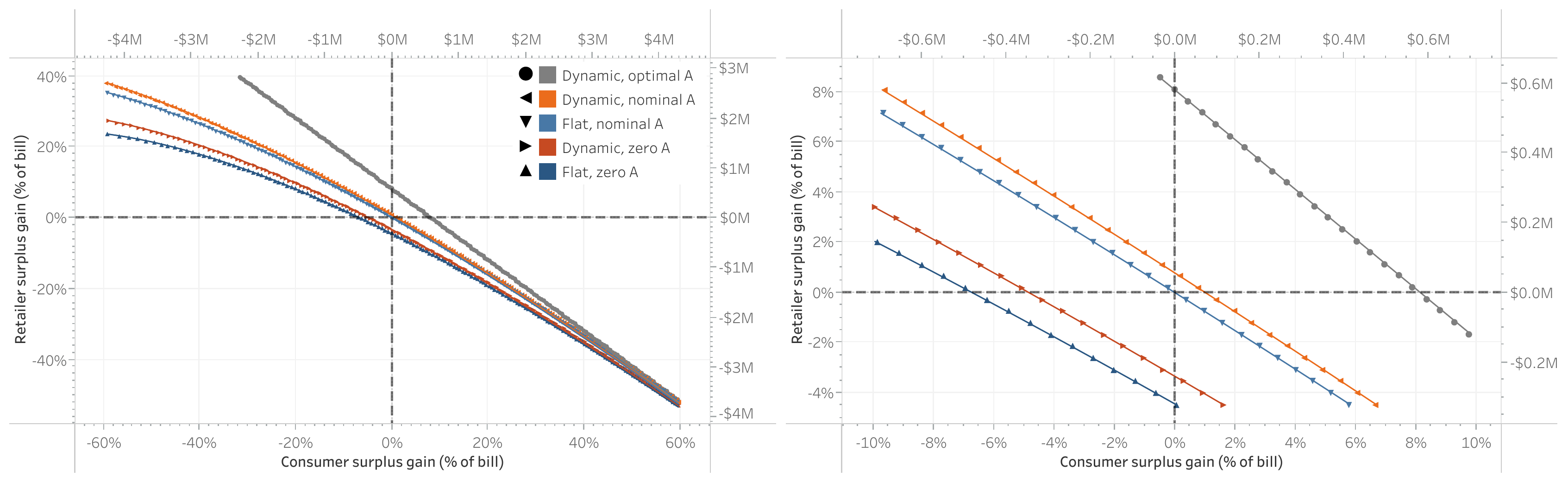}\\[-10pt]
\subcaptionbox{Pareto front for wide range of $F$ including $\overline{\rp}(T^{\textsc{ce}})$.\label{fig:example:Pareto}}
{\rule{0.49\linewidth}{0pt}}
\subcaptionbox{Zoom into neighborhood of $\left(\overline{\cs}(T^{\textsc{ce}}),\overline{\rp}(T^{\textsc{ce}}) \right)$.\label{fig:example:ParetoZoom}}
{\rule{0.49\linewidth}{0pt}}
\caption{CS and RS gains induced by various tariffs parametrized by the RS target $F$ assuming an intermediate own-price elasticity of $\overline{\varepsilon}(\pi^{\textsc{ce}})=-0.3$.}
\label{fig:ParetoFronts:partI}
\end{figure*}

%\begin{figure*}[t]
%\centering
%\fbox{\includegraphics[width=1\linewidth]{\MyLocalPath figures/example/baseCase/ParetoFrontBoth}}\\[-10pt]
%\subfloat[][Pareto front for wide range of $F$ including $\overline{\rp}(T^{\textsc{ce}})$.]{
%\rule{0.5\linewidth}{0pt}\label{fig:example:Pareto}}
%\subfloat[][Zoom into neighborhood of $\left(\overline{\cs}(T^{\textsc{ce}}),\overline{\rp}(T^{\textsc{ce}}) \right)$.]{
%\rule{0.5\linewidth}{0pt}\label{fig:example:ParetoZoom}}
%\caption{CS and RS gains induced by various tariffs parametrized by the RS target $F$ assuming an intermediate own-price elasticity of $\overline{\varepsilon}(\pi^{\textsc{ce}})=-0.3$.}
%\label{fig:ParetoFronts:partI}
%\end{figure*}
	
\else
% Single column version here

\begin{figure*}[t]
\centering
\subfloat[][Pareto front for wide range of $F$ including $\overline{\rp}(T^{\textsc{ce}})$.]{
\includegraphics[width=1\linewidth]{\MyLocalPath figures/example/baseCase/ParetoFront}\label{fig:example:Pareto}}\\
\subfloat[][Zoom into neighborhood of $\left(\overline{\cs}(T^{\textsc{ce}}),\overline{\rp}(T^{\textsc{ce}}) \right)$.]{
\includegraphics[width=1\linewidth]{\MyLocalPath figures/example/baseCase/ParetoFrontZoomIn}\label{fig:example:ParetoZoom}}
\caption{Consumer surplus and retailer surplus gains induced by various tariffs parametrized by the fixed cost parameter $F$ assuming an intermediate own-price elasticity of $\overline{\varepsilon}(\pi^{\textsc{ce}})=-0.3$.}
\label{fig:ParetoFronts}
\end{figure*}

\fi

%% Three figures in line
%\begin{figure*}[t]
%\centering
%\begin{minipage}[b]{\figwid}
%\centering
%{
%\fbox{\includegraphics[height=0.215\paperwidth, trim=0pt 0pt 0pt 0pt, clip=true]{figures/ParetoFront.pdf}}
%}
%\caption{Pareto fronts induced by optimal linear and two-part tariffs.
%The slope $\partial \overline{\rp}(F)/ \partial \overline{\cs}(F) = -\eta$ relates to $\gamma \geq 1$ in \eqref{eq:linearTariff} through $\eta=\frac{2}{1+\gamma}$.}
%\label{fig:noIntegrationGeometry}
%\end{minipage}\hfill
%\begin{minipage}[b]{0.638\linewidth}
%\centering
%\subfloat[][$F \in \left(\overline{\rp}(\pi^*),\overline{\rp}(\pi^\textsc{m}) \right)$.]{
%\label{fig:example:Pareto}
%\fbox{\includegraphics[height=0.2\paperwidth, trim=41pt 0pt 0pt 35pt, clip=false]{figures/ParetoFront-Normalized-noDER-M2200000-alpha2e-001.eps}}
%}\hfill 
%\subfloat[][Zoom into positive quadrant of \subref{fig:example:Pareto}.]{
%\label{fig:example:ParetoZoom}
%\fbox{\includegraphics[height=0.2\paperwidth, trim=41pt 0pt 0pt 35pt, clip=false]{figures/ParetoFront-Normalized-noDER-M2200000-alpha2e-001-ZoomIn.eps}}
%}
%\caption{CS and RS gains induced by various tariffs parametrized by the RS target $F$ assuming an intermediate own-price elasticity of demand of $\varepsilon=-0.3$.}
%\end{minipage}
%%\caption{Pareto fronts induced by deterministic linear and two-part tariffs without renewables or storage integration.}
%\label{fig:ParetoFronts}
%\end{figure*}

We now estimate the performance of
%and \red{distribution}
the optimal linear and two-part day-ahead tariffs in a practical setting.
%Similarly, we study the case where the same levels of DERs are integrated by the retailer.
%We draw conclusions by comparing the efficiency and consumer surplus gains brought by DERs under both types of tariffs and integration models.
Using publicly available data from ConEdison (New York City's largest utility) and NYISO for the 2015 Summer season, we estimate the average daily gains in consumer surplus that both tariffs would have brought relative to the utility's default two-part tariff with flat rate.
Here we assume a linear demand model \eqref{eq:demand function linear}-\eqref{eq:benefit function quadratic} and day-ahead linear and two-part tariffs.

%We use the linear demand response model as follows.
The utility's monthly residential energy sales\footnote{For June through August 2015, which can be found in the EIA-826 database at \url{http://www.eia.gov/electricity/data/eia826/}.} and an estimated residential hourly load profile\footnote{For a residential building in NYC, available in the NREL OpenEI building load database \url{http://en.openei.org/datasets/files/961/pub/}.
The location/model with ID 725033-TMY3-BASE was used.} were used to obtain $2,208$ ($92$ days $\times$ $24$ hr) aggregate hourly consumption data points.
We used these points as iid realizations $\left\{ \b{x}_j \right\}_{j=1}^{92}$ of the random vector $D(\b{1}\pi^{\textsc{ce}},\omega) = \omega - G \b{1} \pi^{\textsc{ce}} \in \Rmbb^{24}$, where $\pi^{\textsc{ce}}$ is ConEdison's flat rate\footnote{NYC's residential default flat rate during Jun-Aug 2015 was $\pi^{\textsc{ce}} = 17.2$ cents/kWh ($47.7\%$ of which correspond to supply charges and $53.3$ to delivery charges).
% ($8.2259$ cents/kWh for supply charges and $9.017$ cents/kWh for delivery charges)
Available at \url{http://www.coned.com/rates/supply_charges.asp} and \url{http://www.coned.com/documents/elecPSC10/SCs.pdf}.}, to fit $G$ and then to estimate $\overline{\omega}$.
Due to the its low-dimensional structure, fitting $G$ can be reduced to determining a scaling parameter after assuming certain own-price elasticity of demand at $\pi^{\textsc{ce}}$\footnote{This is facilitated by assuming a homogeneous thermal parameter $\alpha^i$ across customers, which implies some loss of customer heterogeneity. See \cite{JiaTong16a} for details of the model, where $\alpha^i = 0.2$ is used in a case study.}
$$\overline{\varepsilon}(\pi^{\textsc{ce}}) = \frac{\partial \Embb\left[ D(\b{1}\pi^{\textsc{ce}},\omega)^{\Top} \b{1} \right]/\partial \pi^{\textsc{ce}}}{\Embb\left[D(\b{1}\pi^{\textsc{ce}},\omega)^{\Top} \b{1} \right]/\pi^{\textsc{ce}}} = \frac{\b{1}^{\Top} G \b{1}}{\hat{\b{x}}^{\Top} \b{1} / \pi^{\textsc{ce}}}, $$
where $\hat{\b{x}}$ denotes the sample mean of $\left\{ \b{x}_j \right\}$. %:=\frac{1}{92}\sum_{j=1}^{92} \b{x}_j$.
%As a motivating example, we assume the value $\varepsilon=-0.3$, which is a reasonable estimate of the short-term own-price elasticity of electricity demand \cite{Borenstein05, EPRI08}.
While here we assume a value of $\overline{\varepsilon}(\pi^{\textsc{ce}})=-0.3$, which is a reasonable estimate of the short-term own-price elasticity of electricity demand \cite{EPRI08}, a sensitivity analysis over $\overline{\varepsilon}(\pi^{\textsc{ce}})$ is presented in \cite{Munoz16}. %Lijesen07, Borenstein05
%Because welfare effects induced by tariff changes generally depend on the price elasticity of demand, a sensitivity analysis over $\varepsilon$ was conducted within the range $[-.2,-.6]$.
%According to a survey of empirical estimates \cite{EPRI08}, this is a reasonable range for short-term own-price elasticities.
We further assumed a total of $M=2.2$ million of residential customers and used ConEdison's default residential connection charge, which amounts to $A^{\textsc{ce}}=0.52$ \$/day, to roughly estimate the utility's average daily revenue from residential customers as $\overline{\textsf{rev}}(T^{\textsc{ce}})=\Embb[D(\b{1}\pi^{\textsc{ce}},\omega)]^{\Top} \b{1}\pi^{\textsc{ce}} + MA^{\textsc{ce}}$ or $\$7.19$ million USD.
As for the prices $\lambda$, we used the day-ahead prices for NYC as iid realizations to estimate $\overline{\lambda}$ and $\Sigma_{\lambda,\omega}$ with sample mean and covariance estimators.

We plot in Fig. \ref{fig:example:Pareto} the Pareto fronts
\begin{align} \label{eq:Pareto front plots}
\left\{ \left( \Delta \overline{\cs}(F),\Delta\overline{\rp}(F) \right) \left| \ F \in [\overline{\phi}(\overline{\lambda}), \overline{\phi}(\pi^{\textsc{m}})] \right. \right\}
\end{align}
induced by the optimal linear and two-part tariffs and three other relevant tariffs that satisfy the revenue sufficiency constraint: an optimized linear \emph{flat} tariff with rate $\pi^{\textsc{f}}(F)$, an optimized two-part tariff with \emph{fixed} connection charge $A^{\textsc{ce}}$ and rate $\pi^{\dagger}(F-A^{\textsc{ce}}M)$, and an optimized two-part \emph{flat} tariff with \emph{fixed} connection charge $A^{\textsc{ce}}$ and rate $\pi^{\textsc{ce}}+\Delta(F)$.
The latter can be thought as ConEdison's \emph{adjusted} tariff.
In \eqref{eq:Pareto front plots}, $\Delta \overline{\cs}(F)$ and $\Delta \overline{\rp}(F)$ denote the surplus gains (losses if negative) relative to the corresponding surplus achieved by $T^{\textsc{ce}}$.
For instance, for the optimal two-part tariff $\Delta \overline{\cs}(F) = \overline{\cs}(T^*)-\overline{\cs}(T^{\textsc{ce}})$.
%In Fig. \ref{fig:example:Pareto}, surplus gains are reported as percentages of the (daily) retail surplus associated to the flat rate, $\overline{\rp}(\pi^{flat})$, which is approximately $\$4.69$ million USD for all the values of $\varepsilon$ considered.

Fig. \ref{fig:example:Pareto} compares the tariffs' performances in consumer surplus gains for different retail surplus targets.
At ConEdison's estimated net revenue level $F=\overline{\rp}(T^{\textsc{ce}})$, which corresponds to $\Delta \overline{\rp}(F) = 0$, significant performance differences can be observed among the computed tariffs.
These differences are more evident in Fig. \ref{fig:example:ParetoZoom}, which magnifies Fig. \ref{fig:example:Pareto} around the origin.
%The scale of the axes suggests that these performances differ critically for unrealistically high retail surplus targets, namely more than 1,000 times the base retail surplus, $\overline{\rp}(\pi^{flat})$.
%Conversely, for small retail surplus variations, the gains in CS brought by the optimal linear and two-part tariffs are practically constant with respect to the adjusted flat tariff.
%This is illustrated by Fig. \ref{fig:example:ParetoZoom}, which magnifies Fig. \ref{fig:example:Pareto} and reveals that the Pareto fronts have a slope of nearly -1 for values of $F$ within $\pm1\%$ of $\overline{\rp}(\pi^{flat})$.
%Such slopes imply that the tariff adjustments required to achieve slightly different retail surplus targets induce surplus transfers between consumers and the utility with virtually zero efficiency loses.
%In particular, Fig. \ref{fig:example:ParetoZoom} shows that tariff structure changes do bring noticeable total surplus gains or loses (relative to the adjusted flat tariff) even if the retail surplus target remains unchanged.
Due to its nonzero connection charge, ConEdison's tariff clearly outperforms the tariffs without connection charges, but is outperformed by the other tariffs with connection charges, which are further optimized.
It is particularly interesting that switching to the optimal linear tariff would bring loses in CS ($-4.8\%$ or $-\$345$k USD/day).
Namely, by virtue of a connection charge, even a simple flat tariff can outperform a fairly sophisticated day-ahead hourly volumetric tariff.
Moreover, fully optimizing ConEdison's rate (but not its connection charge) brings rather limited CS gains ($1\%$ or $\$72$k USD/day).
However, switching to the optimal two-part tariff would bring significant gains in CS ($8.1\%$ or $\$582$k USD/day).
This corroborates how effective connection charges can be at increasing the retailer surplus without sacrificing economic efficiency.
%This optimal two-part tariff, which features a connection charge of $A^{*}=2.654$ \$/day, induces no bill changes for a customer $6.62\%$ larger than the average customer and a surplus gain of $0.25$ \$/day.
This optimal two-part tariff, which features a connection charge of $A^{*}=2.65$ \$/day or nearly $80$ \$/month, induces bill reductions for customers $6.62\%$ larger than the average customer and bill increments for all other customers.
Clearly, such charges may be politically unacceptable for low-income customers and may require cross-subsidized reduced tariffs, which have been an industry standard \cite[Sec. 7.4]{BrownSibley86}.
\def\figwid{0.47\linewidth} % For 1 fig per row
\setlength\fboxsep{0pt}
\setlength\fboxrule{0pt}

\section{Conclusions} \label{sec:conclusions}

In this first part, we derive consumer-welfare-maximizing, revenue adequate, and ex-ante linear and two-part dynamic tariffs from the perspective of a regulated retailer.
This initial analysis is for the case without renewables or storage in the distribution system.
Our results generalize previous works by deriving said tariffs for a stochastic and multi-period demand model with intertemporal dependencies and a predetermined lag time between the announcement of the tariff and the beginning of the billing period.
%In particular, we derived the optimal day-ahead two-part tariffs under both integration models and established the extent of their optimality among the larger set of arbitrary (possibly nonlinear) tariffs.
%In particular, we derived the optimal day-ahead two-part tariffs under both integration models and established sufficient conditions for their optimality within the set of arbitrary (nonlinear) tariffs.
%That is, we provided sufficient conditions under which more complex nonlinear tariffs cannot do better than day-ahead two-part tariffs in terms of social welfare.
%The latter result reveals specific features of consumers' behavior that limit the performance of day-ahead tariffs, and thus suggests a type of customers on which the derived tariffs induce a socially optimal behavior.
We established that if the wholesale prices and each customer's consumption are statistically independent, then the optimal two-part tariff is optimal among the class of arbitrary tariffs with the same lag time.

%Our analysis has ignored, for the most part, the consumer welfare distribution effects of the tariffs considered.
While the optimal two-part tariff mitigates inefficiencies induced by the optimal linear tariff, inequity concerns inconsistent with cost causation arise from the structure of the connection charge.
These concerns may become significant with the sparse adoption of behind-the-meter renewables.
%For instance, under such tariff, customers whose consumption is highly correlated with the wholesale prices may be subsidizing those exhibiting milder correlations.
%Our tariff design formulation has ignored common policy objectives beyond efficiency and revenue sufficiency such as fairness (\eg welfare distribution among customers).
%These other criteria are typically in tension with the ones we have considered and make retail tariff design a laborious task (\eg see \cite{JoskowWolfram12, Borenstein11a, Borenstein11b, Borenstein15}).
%Arguably, one of the critical barriers to abrupt retail tariff changes are potential welfare redistribution effects between customers with different income levels \cite{JoskowWolfram12, Borenstein11a, Borenstein11b, Borenstein15}.
While tariff design criteria beyond efficiency and revenue adequacy are out of the scope of our work, it is worth mentioning that allowing discriminatory connection charges can give flexibility to the regulator to achieve different objectives (such as inter-customer cost-causation equity) and provide effective long-term signals (\eg location within the distribution network and investment in on-site generation) \cite{PerezSmeers03}.
%A related example is the active debate on whether utilities should impose or not a higher discriminatory connection charge for customers with behind-the-meter PV solar systems.
%On this regard, our analysis corroborates the idea that for pure efficiency purposes fixed connection charges are preferred to price markups as means to collect nontrivial revenue requirements for the utility; it does not indicate, however, if connection charges should be discriminatory for different types of customers.

%Moreover, our analysis considers only consumption behavior in the short-run in a setting where the level of DER integration is assumed exogenous and its impact on wholesale prices is neglected.
%Our work, however, is a first step towards the more ambitious goal of devising a long-run dynamic model of DER integration derived from first principles.
%%, or to characterize the long-run wholesale market equilibrium in which the level of centralized and decentralized DER becomes endogenous.

%We have ignored other costs that retailers possibly face such as capacity costs.

%Also, when consumer heterogeneity is significant, care needs to be taken when designing two-part tariffs as tradeoffs between efficiency and fairness can emerge as some customers may be better off by choosing not to consume power from the grid at all.
%To mitigate these effects different strategies can be adopted such as offering multiple opt-in tariffs and allowing consumer self-selection \cite[Sec. 5.2]{Wilson93}.

\bibliographystyle{IEEEtran}%IEEEtranN
%\bibliography{Literature}{\markboth{References}{References}}

\begin{thebibliography}{10}
\providecommand{\url}[1]{#1}
\csname url@samestyle\endcsname
\providecommand{\newblock}{\relax}
\providecommand{\bibinfo}[2]{#2}
\providecommand{\BIBentrySTDinterwordspacing}{\spaceskip=0pt\relax}
\providecommand{\BIBentryALTinterwordstretchfactor}{4}
\providecommand{\BIBentryALTinterwordspacing}{\spaceskip=\fontdimen2\font plus
\BIBentryALTinterwordstretchfactor\fontdimen3\font minus
  \fontdimen4\font\relax}
\providecommand{\BIBforeignlanguage}[2]{{%
\expandafter\ifx\csname l@#1\endcsname\relax
\typeout{** WARNING: IEEEtran.bst: No hyphenation pattern has been}%
\typeout{** loaded for the language `#1'. Using the pattern for}%
\typeout{** the default language instead.}%
\else
\language=\csname l@#1\endcsname
\fi
#2}}
\providecommand{\BIBdecl}{\relax}
\BIBdecl

\bibitem{GTM15}
\BIBentryALTinterwordspacing
{Greentech Media (GTM)}, ``Evolution of the grid edge: Pathways to
  transformation,'' White paper, 2015. [Online]. Available:
  \url{https://www.greentechmedia.com}
\BIBentrySTDinterwordspacing

\bibitem{RMI15}
\BIBentryALTinterwordspacing
{Rocky Mountain Institute}, ``The economics of load defection,'' Tech. Rep.,
  2015. [Online]. Available: \url{http://www.rmi.org/}
\BIBentrySTDinterwordspacing

\bibitem{BrownSibley86}
S.~J. Brown and D.~S. Sibley, \emph{The theory of public utility
  pricing}.\hskip 1em plus 0.5em minus 0.4em\relax Cambridge University Press,
  1986.

\bibitem{JoskowWolfram12}
P.~L. Joskow and C.~D. Wolfram, ``Dynamic pricing of electricity,'' \emph{The
  American Economic Review}, vol. 102, no.~3, pp. 381--385, 2012.

\bibitem{BraithwaitEtal07}
S.~Braithwait, D.~Hansen, and M.~O'Sheasy, ``{Retail electricity pricing and
  rate design in evolving markets},'' Edison Electric Institute, Tech. Rep.
  July, 2007.

\bibitem{LBNL16b}
L.~Wood, J.~Howat, R.~Cavanagh, and S.~Borenstein, ``Recovery of utility fixed
  costs: Utility, consumer, environmental and economist perspectives,'' LBNL,
  Tech. Rep., 2016.

\bibitem{Borenstein14}
\BIBentryALTinterwordspacing
S.~Borenstein. (2014) What's so great about fixed charges? Energy Institute at
  HAAS Blog. [Online]. Available: \url{https://energyathaas.wordpress.com}
\BIBentrySTDinterwordspacing

\bibitem{DengEtAl15}
R.~Deng, Z.~Yang, M.-Y. Chow, and J.~Chen, ``{A Survey on Demand Response in
  Smart Grids: Mathematical Models and Approaches},'' \emph{Industrial
  Informatics, IEEE Trans. on}, vol.~11, no.~3, pp. 1--1, 2015.

\bibitem{VardakasEtAl15}
J.~S. Vardakas, N.~Zorba, and C.~V. Verikoukis, ``{A Survey on Demand Response
  Programs in Smart Grids: Pricing Methods and Optimization Algorithms},''
  \emph{Communications Surveys {\&} Tutorials, IEEE}, vol.~17, no.~1, pp.
  152--178, jan 2015.

\bibitem{FaruquiEtal12}
A.~Faruqui, R.~Hledik, and J.~Palmer, ``{Time-Varying and Dynamic Rate
  Design},'' The Brattle Group, Tech. Rep. July, 2012.

\bibitem{BorensteinEtal02}
S.~Borenstein, M.~Jaske, and A.~Rosenfeld, ``Dynamic pricing, advanced
  metering, and demand response in electricity markets,'' Center for the Study
  of Energy Markets, Tech. Rep., 2002.

\bibitem{Borenstein05}
S.~Borenstein, ``{The long-run efficiency of real-time electricity pricing},''
  \emph{Energy Journal}, vol.~26, no.~3, pp. 93--116, 2005.

\bibitem{BorensteinHolland05}
S.~Borenstein and S.~Holland, ``On the efficiency of competitive electricity
  markets with time-invariant retail prices,'' \emph{The Rand Journal of
  Economics}, vol.~36, no.~3, p. 469, 2005.

\bibitem{JoskowTirole06}
P.~Joskow and J.~Tirole, ``Retail electricity competition,'' \emph{The {RAND}
  Journal of Economics}, vol.~37, no.~4, pp. 799--815, 2006.

\bibitem{MohsenianLeon10}
A.~H. Mohsenian-Rad and A.~Leon-Garcia, ``Optimal residential load control with
  price prediction in real-time electricity pricing environments,''
  \emph{{Smart Grids, IEEE Trans. on}}, vol.~1, no.~2, pp. 120--133, 2010.

\bibitem{MohsenianWong10}
A.~H. Mohsenian-Rad, V.~W.~S. Wong, J.~Jatskevich, R.~Schober, and
  A.~Leon-Garcia, ``Autonomous demand-side management based on game-theoretic
  energy consumption scheduling for the future smart grid,'' \emph{{Smart Grid,
  IEEE Trans. on}}, vol.~1, no.~3, pp. 320--331, 2010.

\bibitem{DaryanianEtAl89}
B.~Daryanian, R.~Bohn, and R.~Tabors, ``{Optimal demand-side response to
  electricity spot prices for storage-type customers},'' \emph{{Power Systems,
  IEEE Trans. on}}, vol.~4, no.~3, 1989.

\bibitem{ShaoEtAl11}
S.~Shao, M.~Pipattanasomporn, and S.~Rahman, ``{Demand response as a load
  shaping tool in an intelligent grid with electric vehicles},'' \emph{{Smart
  Grid, IEEE Trans. on}}, vol.~2, no.~4, pp. 624--631, 2011.

\bibitem{MathieuEtAl15}
J.~L. Mathieu, M.~Kamgarpour, J.~Lygeros, G.~Andersson, and D.~S. Callaway,
  ``{Arbitraging intraday wholesale energy market prices with aggregations of
  thermostatic loads},'' \emph{{Power Systems, IEEE Trans. on}}, vol.~30,
  no.~2, pp. 763--772, 2015.

\bibitem{LiEtAl16}
S.~Li, W.~Zhang, J.~Lian, and K.~Kalsi, ``Market-based coordination of
  thermostatically controlled loads--{Part I: A} mechanism design
  formulation,'' \emph{{Power Systems, IEEE Trans. on}}, vol.~31, no.~2, 2016.

\bibitem{ChenEtal12}
L.~Chen, N.~Li, L.~Jiang, and S.~H. Low, ``Optimal demand response: problem
  formulation and deterministic case,'' in \emph{Control and optimization
  methods for electric smart grids}.\hskip 1em plus 0.5em minus 0.4em\relax
  Springer, 2012, ch.~3, pp. 63--85.

\bibitem{LiEtAl13}
C.~Li, S.~Tang, Y.~Cao, Y.~Xu, Y.~Li, J.~Li, and R.~Zhang, ``{A New Stepwise
  Power Tariff Model and Its Application for Residential Consumers in Regulated
  Electricity Markets},'' \emph{{Power Systems, IEEE Trans. on}}, vol.~28,
  no.~1, pp. 300--308, 2013.

\bibitem{YangEtAl13}
P.~Yang, G.~Tang, and A.~Nehorai, ``{A game-theoretic approach for optimal
  time-of-use electricity pricing},'' \emph{{Power Systems, IEEE Trans. on}},
  vol.~28, no.~2, pp. 884--892, 2013.

\bibitem{MengZeng13}
F.-L. Meng and X.-J. Zeng, ``{A Stackelberg game-theoretic approach to optimal
  real-time pricing for the smart grid},'' \emph{Soft Computing}, vol.~17,
  no.~12, pp. 2365--2380, 2013.

\bibitem{BinaAhmadi15}
M.~T. Bina and D.~Ahmadi, ``{Stochastic Modeling for the Next Day Domestic
  Demand Response Applications},'' \emph{{Power Systems, IEEE Trans. on}},
  vol.~30, no.~6, pp. 2880--2893, 2015.

\bibitem{JiaTong16a}
L.~Jia and L.~Tong, ``Dynamic pricing and distributed energy management for
  demand response,'' \emph{{Smart Grid, IEEE Trans. on}}, vol.~7, 2016, earlier
  arXiv version: \url{http://arxiv.org/abs/1601.02319}.

\bibitem{MunozTong16}
\BIBentryALTinterwordspacing
D.~Munoz-Alvarez and L.~Tong, ``Distributed renewables and storage and the
  efficiency of connection charges,'' in \emph{Information Theory and
  Applications Workshop (ITA)}, 2016, pp. 1--6. [Online]. Available:
  \url{http://ita.ucsd.edu/workshop/16/files/paper/paper_4289.pdf}
\BIBentrySTDinterwordspacing

\bibitem{CrewEtAl95}
M.~A. Crew, C.~S. Fernando, and P.~R. Kleindorfer, ``{The theory of peak-load
  pricing: A survey},'' \emph{Journal of Regulatory Economics}, vol.~8, no.~3,
  pp. 215--248, nov 1995.

\bibitem{AguirreEtal10}
I.~Aguirre, S.~Cowan, and J.~Vickers, ``{Monopoly price discrimination and
  demand curvature},'' \emph{American Economic Review}, vol. 100, no.~4, pp.
  1601--1615, 2010.

\bibitem{Burgeretal07}
M.~Burger, B.~Graeber, and G.~Schindlmayr, \emph{Managing Energy Risk: An
  Integrated View on Power and Other Energy Markets}, ser. The Wiley Finance
  Series.\hskip 1em plus 0.5em minus 0.4em\relax John Wiley \& Sons, 2007, vol.
  425.

\bibitem{PerezSmeers03}
I.~J. P\'erez-Arriaga and Y.~Smeers, ``{Guidelines on Tariff Setting},'' in
  \emph{Transport Pricing of Electricity Networks}.\hskip 1em plus 0.5em minus
  0.4em\relax Springer US, 2003, no.~l, pp. 175--203.

\bibitem{MunozTong16partIIarxiv}
D.~Munoz-Alvarez and L.~Tong, ``{On the Efficiency of Dynamic Retail Tariffs
  with Connection Charges---Part II: Distributed Renewable and Storage
  Resources},'' \emph{ArXiv Preprint}, 2017.

\bibitem{EPRI08}
{Electric Power Research Institute (EPRI)}, ``{Price Elasticity of Demand for
  Electricity : A Primer and Synthesis},'' Tech. Rep., 2008.

\bibitem{Munoz16}
D.~{Mu\~noz-\'Alvarez}, ``Regulatory approaches to the integration of renewable
  and storage resources into electricity markets,'' Ph.D. dissertation, Cornell
  University, 2017.

\end{thebibliography}

\appendix
%\appendix
\label{sec:appendix:part1}

%\subsection{Remarks on Consumer Model} \label{app:consumer model}
%\blue{
%The optimality condition assumed to relate $D^i$ and $S^i$, or better, to characterize the demand function $D^i$ using $S^i$ as a more primitive characterization of each customer can be rigorously justified as follows...
%}
%
%
%\propositionAlt{ \label{prop:FOC demand function}
%The demand function $D^i(\pi,\omega^i)$ satisfies the first-order necessary condition
%\begin{align} \label{eq:FOC consumer problem}
%\Embb\left[ v^i_k(D^i(\pi,b^i),b^i) \ | \ \omega^i_1,\ldots,\omega^i_k \right] = \pi_k
%\end{align}
%%$$ \Embb_{|k}[v^i_k(q^i(T,b^i),b^i)] = \Embb_{|k}[\partial T(q^i(T,b^i))/ \partial q^i_k,$$
%for $k=1,\ldots,N$, where the expectation is taken over $\omega^i$ conditioned on $(\omega^i_1,\ldots,\omega^i_k)$ and $v^i_k = \partial S^i / \partial q^i_k$ is the marginal benefit of consuming an incremental unit at time $k$.
%}

\subsubsection*{Proof of Theorem \ref{thm:affineTariff}} \label{proof: Thm 1}

Solving $\overline{\rp}(T)=F$ in \eqref{eq:reg problem} for $A$ and substituting in the objective $\overline{\cs}(T)$ yields the expression
$$
\overline{\sw}(T) \ = \ \sum_{i=1}^{M} \Embb\left[ S^i(D^{i}(\pi,\omega^i), \omega^i) -  \lambda^{\Top} D^{i}(\pi,\omega^i) \right].
$$
Differentiating the objective $\overline{\sw}(T)$ over $\pi$ yields 
\begin{align}
\nabla_{\pi} \overline{\sw}(T) &= \sum_{i=1}^{M} \Embb\left[ \nabla_{\pi}D^i(\pi,\omega^i) \left( \nabla_{q}S^i(D^{i}(\pi,\omega^i), \omega^i) - \lambda \right) \right] \nn\\
	&= \Embb\left[ \nabla_{\pi}D(\pi,\omega) ( \pi - \lambda ) \right] \label{eq:proofThm1}
\end{align}
where the second equality follows from Prop. \ref{prop:implication of FOC demand function} below.
Now, the assumption that each $\nabla_{\pi}D^i(\pi,\omega^i)$ is negative definite implies that $\Embb [ \nabla_{\pi}D(\pi,\omega)]$ is invertible.
Hence, one can rewrite the necessary first-order condition (FOC), which is obtained from equating $ \nabla_{\pi} \overline{\sw}(T) = 0$, as \eqref{eq:affineTariff}.
%Now, Assumption \ref{Assumption 1} readily implies that the Jacobian of \eqref{eq:proofThm1}, say $\nabla^2 \overline{\sw}(\pi) \in \Rmbb^{N\times N}$, is negative semidefinite, thus implying the concavity of $\overline{\sw}(\pi)$ in $\pi$.
Given the expression \eqref{eq:affineTariff} for $\pi^*$, \eqref{eq:connection charge} can be obtained by solving the equality constraint $\overline{\rp}(T^*)=F$ for $A^*$.
Lastly, these FOC are sufficient for the optimality of $(A^*,\pi^*)$ since $\overline{\sw}(T) = \overline{\cs}(T) + \overline{\rp}(T)$ is strictly concave in $\pi$ according to Prop. \ref{prop:concavity rs}.
\hfill $\blacksquare$

\propositionAlt{ \label{prop:FOC demand function}
For each $k=1,\ldots,N$, $D^i(\cdot,\cdot)$ satisfies
%$
%\mbox{$\Embb\left[\left. \frac{\partial S^i(D^i(\pi,\omega^i),\omega^i)}{\partial q^i_k} \right| \omega^i_1,\ldots,\omega^i_k \right]$} = \pi_k
%$
\begin{align} \label{eq:FOC consumer problem}
\Embb\left[\left. \partial S^i(D^i(\pi,\omega^i),\omega^i)/\partial q^i_k \right| \ \omega^i_1,\ldots,\omega^i_k \right] = \pi_k
\end{align}
where the conditional expectation is taken over $\omega^i$ conditioned on $\omega^i_1,\ldots,\omega^i_k$.
}
\proof{
\eqref{eq:FOC consumer problem} are FOCs of customer $i$'s multi-stage decision problem \eqref{eq:consumer surplus individual} of sequentially choosing $q^i_1,\ldots,q^i_k,\ldots,q^i_N$ contingent on the so-far observed local states $\omega^i_1,\ldots,\omega^i_k$ to maximize his expected net utility or surplus, \ie
\begin{align} \label{eq:consumer surplus individual 2}
\max _{q^i(\cdot)} &\ \Embb_{\omega^i} \left[ S^i \left( q^i(\omega^i), \omega^i \right) - T\left( q^i(\omega^i) \right) \right]
\end{align}
given a tariff $T(q^i)=A+\pi^{\Top}q^i$ known in advance (ex-ante).
Hence, the optimal solution $q^i(\cdot) = D^i(\pi,\cdot)$ of \eqref{eq:consumer surplus individual 2} must satisfy these necessary optimality conditions.

Each of these stationarity (KKT) FOCs is obtained by differentiating the objective in \eqref{eq:consumer surplus individual 2} with respect to $q^i_k$ and equating the result to zero.
To see that, note that for each time $k=1,\ldots,N$ one can use the law of total expectation to rewrite the objective in terms of an expectation conditioned on the information observed up to $k$, \ie
\begin{align*}
\Embb_{\omega^i_1,\ldots,\omega^i_k} \left[
\Embb \left[
\left.
S^i \left( q^i(\omega^i), \omega^i \right) - T\left( q^i(\omega^i) \right)
\right| \omega^i_1,\ldots,\omega^i_k \right]
\right].
\end{align*}
The FOC in \eqref{eq:FOC consumer problem} follows since $\partial T(q^i(\omega^i))/\partial q^i_k = \pi_k$.
\hfill $\blacksquare$
}

\propositionAlt{ \label{prop:implication of FOC demand function}
$D^i(\cdot,\cdot)$ satisfies the equation
\begin{align} \label{eq:prop implication of FOC}
\partial \Embb[S^i(D^i(\pi,\omega^i),\omega^i)]/\partial \pi
&= \Embb\left[ \nabla_{\pi}D^i(\pi,\omega^i) \right] \pi
\end{align}
where the expectations are taken over $\omega^i$.
}
\proof{
We establish this vectorial identity proving each component separately.
Assuming the differentiation and expectation operators can be exchanged, we apply the chain rule to the $k$-th component of the left-hand-side of \eqref{eq:prop implication of FOC} yielding the following sequence of equalities
\begin{align*}
\frac{\partial \Embb[S^i(D^i(\pi,\omega^i),\omega^i)]}{\partial \pi_k} 
&=
\Embb\left[ \sum_{t=1}^N
\frac{\partial S^i}{\partial q^i_t} \cdot
\frac{\partial D^i_t}{\partial \pi_k}
\right] \\
&=
\sum_{t=1}^N
\Embb\left[
\Embb \left[ \left.
\frac{\partial S^i}{\partial q^i_t} \cdot
\frac{\partial D^i_t}{\partial \pi_k}
\right| \omega^i_1,\ldots,\omega^i_t \right]
\right] \displaybreak[0] \\
&=
\sum_{t=1}^N
\Embb\left[
\Embb \left[ \left.
\frac{\partial S^i}{\partial q^i_t}
\right| \omega^i_1,\ldots,\omega^i_t \right]
\frac{\partial D^i_t}{\partial \pi_k}
\right] \displaybreak[0]  \\
&=
\pi^{\Top}
\Embb\left[
\frac{\partial D^i}{\partial \pi_k}
\right],
\end{align*}
where the second equality follows from the law of total expectation, the third equality is due to the causality of $D^i(\pi,\omega^i)$, and the last equality is due to Prop. \ref{prop:FOC demand function}.
The identity \eqref{eq:prop implication of FOC} readily follows.

% Strict convexity of cs in \pi (Uses two other propositions)
\propositionAlt{ \label{prop:convexity cs}
For any affine tariff $T(q) = A+\pi^{\Top} q$, $\overline{\cs}(T)$ is strictly convex and (componentwise) decreasing in $\pi$.
}
\proof{
Differentiating $\overline{\cs}(T)$ with respect to $\pi$ (using the chain rule) and leveraging on the expression \eqref{eq:prop implication of FOC} from Prop. \ref{prop:implication of FOC demand function} one readily obtains
$
\nabla_{\pi} \overline{\cs}(T)	= -\Embb[ D(\pi,\omega)].
$
Because $D^i(\pi,\omega^i)$ is nonnegative, the expression above implies that $\overline{\cs}(T)$ is componentwise decreasing $\pi$.
Moreover, we have that
\begin{align} \label{eq:cs Hessian}
\nabla^2_{\pi} \overline{\cs}(T)	= -\Embb[ \nabla_{\pi} D(\pi,\omega)].
\end{align}
Recall that a function is strictly convex (over a convex domain) if its Hessian is positive definite (over said domain).
Hence, the strict convexity of $\overline{\cs}(T)$ in $\pi$ readily follows from the assumed negative definiteness of each matrix $\nabla_{\pi} D^i(\pi,\omega)$.
\hfill $\blacksquare$
}

% Strict concavity of weigthed surplus in \pi
\propositionAlt{ \label{prop:concavity rs}
Consider an affine tariff $T(q) = A+\pi^{\Top} q$.
If Assumption \ref{Assumption 1} holds then
$\overline{\rp}(T)$ and 
the weighted surplus $\overline{\cs}(T)+ \gamma \overline{\rp}(T)$ are strictly concave in $\pi$ for any $\gamma \geq 1$.
}
\proof{
First, we differentiate twice $\overline{\rp}(T)$ with respect to $\pi$.
Using Prop. \ref{prop:implication of FOC demand function} and some algebra one obtains the Hessian
\begin{align} \label{eq:rs Hessian}
\nabla^2_{\pi} \overline{\rp}(T) = \Embb[\nabla_{\pi} D(\pi,\omega)] + \nabla g(\pi).
\end{align}
Recall that a function is strictly concave (over a convex domain) if its Hessian is negative definite (over said domain).
Hence, the assumed negative definiteness of $\nabla_{\pi} D^i(\pi,\omega)$ and of $\nabla g(\pi)$ (Assumption \ref{Assumption 1}) readily imply the negative definiteness of $\nabla^2_{\pi} \overline{\rp}(T)$ and thus the strict concavity of $\overline{\rp}(T)$ in $\pi$.

%Recall now that a function is concave over a convex domain if and only if its Hessian is negative semidefinite (nsd) over said domain.
%Hence, the assumed negative definiteness of $\nabla_{\pi} D(\pi,\omega)$ and negative semidefiniteness of $\nabla g(\pi)$ (Assumption \ref{Assumption 1}) imply the negative semidefiniteness of $\nabla^2_{\pi} \overline{\rp}(T)$ and thus the concavity of $\overline{\rp}(T)$ in $\pi$.

It remains to show that $\nabla^2_{\pi} (\overline{\cs}(T)+ \gamma \overline{\rp}(T))$ is negative definite for $\gamma \geq 1$.
From \eqref{eq:cs Hessian} and \eqref{eq:rs Hessian} we have that
\begin{align*}
\nabla^2_{\pi} (\overline{\cs}(T)+ \gamma \overline{\rp}(T)) = (\gamma-1) \Embb[\nabla_{\pi}D(\pi,\omega)] + \gamma \nabla g(\pi).
\end{align*}
Similar arguments yield the desired result since both terms on the right hand side are negative definite for $\gamma \geq 1$.
\hfill $\blacksquare$
}

\subsubsection*{Proof of Corollary \ref{cor:affineTariff1}} \label{proof: Cor 1}

If $\nabla_{\pi} D(\pi,\omega)$ and $\lambda$ are uncorrelated, then clearly $\Embb[\nabla_{\pi} D(\pi,\omega) \lambda] = \Embb[\nabla_{\pi} D(\pi,\omega)] \overline{\lambda}$.
It follows that the expression for $\pi^*$ in \eqref{eq:affineTariff} reduces to $\pi^*=\overline{\lambda}$.
In turn, the expression for $\overline{\phi}(\pi^*)$ in \eqref{eq:connection charge} reduces to 
$
\overline{\phi}(\pi^*) = \overline{\phi}(\overline{\lambda}) = -\tr(\cov(\lambda,D(\overline{\lambda},\omega))),
$
thus readily simplifying \eqref{eq:connection charge} to
\begin{align*}
A^* = \mbox{$\frac{1}{M}$} \left( F + \tr \left( \cov \left( \lambda,D \left( \overline{\lambda},\omega \right) \right) \right) \right).
\qedadhoc
\end{align*}

\subsubsection{Proof of Corollary \ref{cor:affineTariff2}} \label{proof: Cor 2}

Leveraging on \eqref{eq:cs} and \eqref{eq:rp}, the expected total surplus induced by the tariff $T^*$ characterized in Thm. \ref{thm:affineTariff} can be written as
\begin{align*}
\overline{\sw}^*	&= \overline{\sw}(T^*) \\
					&= \overline{\cs}(T^*) + \overline{\rp}(T^*) \\
					&= \sum_{i=1}^{M} \Embb\left[ S^i(D^{i}(\pi^*,\omega^i), \omega^i) -  \lambda^{\Top} D^{i}(\pi^*,\omega^i) \right].
\end{align*}
Clearly, $\overline{\sw}(T)$ is a function of $\pi$ but not of $F$ or $A$.
It remains to show that $\pi^*$ does not depend on $F$.
But the FOC characterizing $\pi^*$ is obtained by differentiating $\overline{\sw}(T)$ with respect to $\pi$ and equating it to zero.
It follows that $\overline{\sw}^*$ does not on $F$.
Moreover, since $\overline{\rp}(T^*)=F$ must hold at optimality, it readily follows that
\begin{align*}
\overline{\cs}(T^*) = \overline{\sw}(T^*) - \overline{\rp}(T^*) = \overline{\sw}^* - F.
\qedadhoc
\end{align*}

%Theorem \ref{thm:affineTariff} implies that, unlike $A^*$, $\pi^*$ is independent of $F$, and so is the expected total surplus $\overline{\sw}^*:=\overline{\sw}(T^*) = \overline{\cs}(T^*)+\overline{\rp}(T^*)$ since the terms depending on $A^*$ in $\overline{\cs}(T^*)$ and $\overline{\rp}(T^*)$ cancel each other.
%Since $\overline{\rp}(T^*)=F$ at optimality, it follows that $\overline{\cs}^*(F):= \overline{\cs}(T^*)= \overline{\sw}^* - F$.

\subsubsection{Proof of Theorem \ref{thm:linearTariff}} \label{proof: Thm 2}

Consider the Lagrangian of problem \eqref{eq:reg problem}
\begin{align} \label{eq:Lagrangian}
\Lmsc (\pi,\gamma) = \overline{\cs}(T) + \gamma (\overline{\rp}(T)-F),
\end{align}
where $\gamma \in \Rmbb$ is the multiplier of the equality constraint.
Differentiating $\Lmsc (\pi,\gamma)$ over $\pi$ yields
\begin{align}
\nabla_{\pi} \Lmsc &= \nabla_{\pi}\overline{\cs}(T) + \gamma \nabla_{\pi}\overline{\rp}(T) \nn\\
	&= \Embb[ \gamma \nabla_{\pi} D(\pi,\omega)(\pi-\lambda) + (\gamma-1)D(\pi,\omega) ], \label{eq:proof Thm 2}
\end{align}
where we use Prop. \ref{prop:implication of FOC demand function} to obtain the second equality.
The necessary FOC in \eqref{eq:linearTariff} follows from equating \eqref{eq:proof Thm 2} to zero after some algebra.
These operations use the negative definiteness and thus invertibility of the matrix $\Embb[ \nabla_{\pi} D(\pi,\omega)]$ and the expression for $\pi^*$ in \eqref{eq:affineTariff} (Thm. \ref{thm:affineTariff}).
%The equivalent FOC in \eqref{eq:linearTariffb} is a component-wise transformation of \eqref{eq:linearTariff}.

The expressions in \eqref{eq:linearTariffb} characterizing $\pi^{\dagger}$ componentwise can be obtained from \eqref{eq:linearTariff} after some algebraic manipulations using the definition of price elasticity in \eqref{eq:elasticity}.
In particular, subtracting $\pi^*$ from both sides of \eqref{eq:linearTariff} and left-multiplying by $\Embb[ \nabla_{\pi} D(\pi^{\dagger},\omega)]$ one obtains
\begin{align*}
%\pi^{\dagger} &= \pi^{*} -  \mbox{$\frac{\gamma-1}{\gamma}$} \Embb[\nabla_{\pi} D(\pi^{\dagger},\omega)]^{-1} \Embb[ D(\pi^{\dagger},\omega)] \\
%\pi^{\dagger} - \pi^* &= -  \mbox{$\frac{\gamma-1}{\gamma}$} \Embb[\nabla_{\pi} D(\pi^{\dagger},\omega)]^{-1} \Embb[ D(\pi^{\dagger},\omega)] \\
\Embb[\nabla_{\pi} D(\pi^{\dagger},\omega)] (\pi^{\dagger} - \pi^*) &= -  \mbox{$\frac{\gamma-1}{\gamma}$} \Embb[ D(\pi^{\dagger},\omega)].
\end{align*}
Componentwise, for each $k=1,\ldots,N$, we have
\begin{align*}
\Embb[\nabla_{\pi} D_k(\pi^{\dagger},\omega)]^{\Top} (\pi^{\dagger} - \pi^*) &= -  \mbox{$\frac{\gamma-1}{\gamma}$} \Embb[ D_k(\pi^{\dagger},\omega)] \\
\frac{\sum_{t=1}^N \Embb[\partial D_k(\pi^{\dagger},\omega)/\partial \pi_t] (\pi^{\dagger}_t - \pi^*_t)}{\Embb[ D_k(\pi^{\dagger},\omega)]} &= -\frac{\gamma-1}{\gamma} \\
\sum_{t=1}^{N} -\overline{\varepsilon}_{kt}(\pi^{\dagger}) \left( \frac{\pi^{\dagger}_t - \pi^*_t}{\pi^{\dagger}_t} \right)  &= \frac{\gamma-1}{\gamma}.
\end{align*}

%This is because Assumption \ref{Assumption 1} and the negative definiteness of $\nabla_{\pi} D(\pi,\omega)$ imply that $\nabla_{\pi} \Lmsc$ is decreasing in $\pi$ for $\gamma \geq 1$ and increasing for $\gamma \leq 0$.
%Hence, $\Lmsc (\pi,\gamma)$ is concave in $\pi$ for $\gamma \geq 1$, convex if $\gamma \leq 0$, and neither if $\gamma \in (0,1)$.
%Moreover, when $\gamma=1$, $\pi^{\dagger}=\pi^*$ which implies that $\overline{\sw}(\pi^{\dagger}) = \overline{\sw}(\pi^*) = \overline{\sw}^*$ and $F=\overline{\rp}(\pi^{\dagger}) = \overline{\rp}(\pi^*) = \phi(\pi^*)$.

As for the last statement of the theorem, consider that maximizing $\overline{\rp}(T)$ (or equivalently $\overline{\phi}(\pi)$) over $\pi$ yields the related stationarity FOC
\begin{align} \label{eq:linearTariff3}
\pi^{\textsc{m}} = \pi^{*} -  \Embb[\nabla_{\pi} D(\pi^{\textsc{m}},\omega)]^{-1} \Embb[ D(\pi^{\textsc{m}},\omega)].
\end{align}
This condition, which is similar to \eqref{eq:linearTariff}, characterizes the unregulated monopoly price $\pi^{\textsc{m}}$.
Indeed, \eqref{eq:linearTariff3} can be obtained from \eqref{eq:linearTariff} by replacing $\mbox{$\frac{\gamma-1}{\gamma}$}$ by $1$, or equivalently, by letting $\gamma \rightarrow \infty$ on both sides of \eqref{eq:linearTariff}.

Now, because $\pi^{\textsc{m}}$ maximizes $\overline{\phi}(\pi)$ and $\overline{\rp}(T)$ over $\pi \geq 0$, there does not exist $\pi \geq 0$ such that $\overline{\rp}(T) > \overline{\phi}(\pi^{\textsc{m}})$.
Hence, when restricted to linear tariffs, problem \eqref{eq:reg problem} is infeasible for all $F > \overline{\phi}(\pi^{\textsc{m}})$.
For all the other values of $F$ within the considered regime, \ie $F \in [\overline{\phi}(\pi^*),\overline{\phi}(\pi^{\textsc{m}})]$, the fact that $\pi^*$ and $\pi^{\textsc{m}}$ achieve $\overline{\phi}(\pi^*)$ and $\overline{\phi}(\pi^{\textsc{m}})$, respectively, and the concavity of $\overline{\rp}(T)$ (Prop. \ref{prop:concavity rs}) imply the feasibility of problem \eqref{eq:reg problem} over said interval of $F$.
Concluding, in the regime $F \geq \overline{\phi}(\pi^*)$, problem \eqref{eq:reg problem} over linear tariffs is feasible if $F \leq \overline{\phi}(\pi^{\textsc{m}})$ and unfeasible otherwise.
%Indeed, one can show that the optimal multiplier $\gamma(F) \in [1,\infty]$ is a strictly monotone function of $F \in [\overline{\phi}(\pi^*),\overline{\phi}(\pi^{\textsc{m}})]$.

\subsubsection{Proof of Corollary \ref{cor:linearTariff}} \label{proof: Cor 3}

At optimality, we have that $\overline{\rp}(T^{\dagger}) = F$ and thus $d\overline{\rp}(T^{\dagger})/dF = 1$.
The envelope theorem further implies that the derivative of the value function $\overline{\cs}(T^{\dagger})$ of problem \eqref{eq:reg problem} with respect to the parameter $F$ can be computed as
$
d\overline{\cs}(T^{\dagger})/dF = \partial \Lmsc(\pi^{\dagger},\gamma)/\partial F = -\gamma
$.
%$\mbox{$\frac{d\overline{\cs}^{\dagger}(F)}{dF}$} = -\gamma$.
The total derivative of the expected total surplus with respect to $F$ is then
$d\overline{\sw}(T^{\dagger})/dF = 1-\gamma$.
%$\mbox{$\frac{d\overline{\sw}^{\dagger}(F)}{dF}$} = 1-\gamma$, 
Hence, to show that $\overline{\cs}(T^{\dagger})$ and $\overline{\sw}(T^{\dagger})$ are decreasing and concave functions of $F$ over $F\in[\overline{\phi}(\pi^*),\overline{\phi}(\pi^{\textsc{m}})]$, it suffices to show that $\gamma$, as a function of $F$, satisfies $\gamma \geq 1$ and $d\gamma/dF \geq 0$ over $F\in[\overline{\phi}(\pi^*),\overline{\phi}(\pi^{\textsc{m}})]$.
Moreover, it is clear from the Lagreangian $\Lmsc(\pi,\gamma)$ in \eqref{eq:Lagrangian} that $\overline{\rp}(T^{\dagger})$ must be a decreasing function of the parameter $\gamma$.
Conversely, the constraint $\overline{\rp}(T) = F$ implies that $\overline{\rp}(T^{\dagger})$ is a strictly increasing function of $F$.
Hence, we have that $\gamma$ increases as $F$ increases, and thus $d\gamma/dF \geq 0$.

Finally, note that $\gamma = 1$ and $\pi^{\dagger} = \pi^*$ are optimal for $F = \overline{\phi}(\pi^*)$ since they satisfy the FOC \eqref{eq:linearTariff} and the constraint $\overline{\rp}(T^{\dagger}) = F$.
Recall also from Thm. \ref{thm:linearTariff} that the problem at hand remains feasible if $F \leq \overline{\phi}(\pi^{\textsc{m}})$.
Because $d\gamma/dF \geq 0$, we can conclude that $\gamma \geq 1$ for $F \in [\overline{\phi}(\pi^*) , \overline{\phi}(\pi^{\textsc{m}})]$, thus completing the proof.

\subsubsection{Proof of Theorem \ref{thm:optimality}}
\label{proof:thm:optimality}

We prove this result by showing that the optimal two-part tariff $T^*$ characterized by Theorem \ref{thm:affineTariff} attains an upper bound for the performance of all \emph{ex-ante} tariffs.
This upper bound is the performance achieved by a social planner who directly makes all decisions on behalf of consumers.
The social planner is unlike the regulator who is limited to coordinate such decisions indirectly through a tariff.
To obtain a tight upper bound for ex-ante tariffs only (rather than a looser bound for all possibly ex-post tariffs), the social planner makes customers' decisions relying only on the information observable by each of them (\ie $\omega^i$) as opposed to based on global information (\eg $\xi=(\lambda,\omega^1,\ldots,\omega^M)$).

Consider the social planner's problem in \eqref{eq:social planner}, which corresponds to the regulator's problem \eqref{eq:reg problem} in the absence of any DERs.
Therein, the notation $q^i(\omega^i)$ indicates the restriction of the social planner to make (causal) decisions \emph{contingent only} on the local state of each customer $\omega^i$.
%A discussion of problem \eqref{eq:social planner with DERs} in the absence of DERs can be found in \cite[Sec. III.C]{MunozTong16partI}; notably, it implicitly restricts the social planner's consumption and storage operation decisions to be contingent only on the local state of each customer.
Recall from \eqref{eq:consumer surplus individual} that the expected consumer surplus for a given ex-ante tariff is given by $\overline{\cs}(T) = \sum_{i=1}^M \overline{\cs}^i(T)$, where
\begin{align*}
\overline{\cs}^i(T) &= \max_{q^i(\cdot)} \Embb \left[ S^i\left(q^i(\omega^i),\omega^i\right) - T\left(q^i(\omega^i)\right) \right] \\
&= \Embb \left[ S^i(q^{i}(T,\omega^i),\omega^i) - T\left( q^{i}(T,\omega^i) \right) \right],
%&= \sum_{i=1}^M \Embb \left[ S^i(q^{i}(T,\omega^i),\omega^i) - T(q^{i}(T,\omega^i)) \right]
\end{align*}
and the corresponding expected retailer surplus is given by
\begin{align*}
\overline{\rp}(T) = \Embb \left[ \sum_{i=1}^M T\left(q^{i}(T,\omega^i) \right) - \lambda^{\Top} q^{i}(T,\omega^i) \right],
\end{align*}
and the expected total surplus by 
\begin{align*}
\overline{\sw}(T) &= \overline{\cs}(T)+\overline{\rp}(T) \\
&= \Embb \left[ \sum_{i=1}^M S^i\left(q^{i}(T,\omega^i),\omega^i\right) - \lambda^{\Top} q^{i}(T,\omega^i) \right].
\end{align*}
The following sequence of equalities/inequalities shows that problem \eqref{eq:social planner} provides an upper bound to problem \eqref{eq:reg problem}.
\begin{align} 
 \max_{T(\cdot)} & \{ \overline{\cs}(T) \ | \ \overline{\rp}(T) = F \} + F \nn\\ 
&= \max_{T(\cdot)} \{ \overline{\cs}(T) + F \ | \ \overline{\rp}(T) = F \} \nn\\ 
&= \max_{T(\cdot)} \{ \overline{\cs}(T) + \overline{\rp}(T) \ | \ \overline{\rp}(T) = F \} \nn\\ 
&= \max_{T(\cdot)} \{ \overline{\sw}(T) \ | \ \overline{\rp}(T) = F \} \nn \displaybreak[0] \\ 
&\leq \max_{T(\cdot)} \ \overline{\sw}(T) \label{eq:reg problem UB 1} \\
&\leq \max_{\{q^i(\cdot)\}_{i=1}^M} \Embb\left[ \sum_{i=1}^MS^i(q^i(\omega^i),\omega^i) - \lambda^{\Top}q^i(\omega^i) \right] \label{eq:reg problem UB 2} \\
&= \max_{\{q^i(\cdot)\}_{i=1}^{M}} \ \overline{\sw}. \label{eq:reg problem UB 3}
\end{align}

In particular, the inequality in \eqref{eq:reg problem UB 2} holds because $\overline{\sw}(T)$ depends on $T$ only through $q^{i}(T,\omega^i)$.
This implies that maximizing $\overline{\sw}(T)$ directly over $\{q^i(\cdot)\}_{i=1}^M$ rather than indirectly over $T(\cdot)$ is a relaxation of the optimization in \eqref{eq:reg problem UB 1}.
Clearly, the problem in \eqref{eq:reg problem UB 2} corresponds to the social planner's problem in \eqref{eq:reg problem UB 3} and \eqref{eq:social planner}.

It suffices to show now that $T^*$ attains the upper bound in \eqref{eq:reg problem UB 3}.
To that end, we use the independence sufficient condition $\omega \perp \lambda$.
We show that, under said condition, the expected total surplus $\overline{\sw}(T^*)$ matches the upper bound.
First note that the condition $\omega \perp \lambda$ allows to rewrite the upper bound in \eqref{eq:reg problem UB 2} and \eqref{eq:reg problem UB 3} as follows.
\begin{align*}
\max_{\{q^i(\cdot)\}_{i=1}^{M}} & \Embb_{\xi}\left[ \sw \right] \\
&= \sum_{i=1}^M \max_{q^i(\cdot)} \ \Embb_{\omega^i} \left[ S^i(q^i(\omega^i),\omega^i) - \Embb\left[\lambda^{\Top} | \omega^i \right] q^i(\omega^i) \right] \\
&= \sum_{i=1}^M \max_{q^i(\cdot)} \ \Embb_{\omega^i}\left[ S^i(q^i(\omega^i),\omega^i) - \overline{\lambda}^{\Top} q^i(\omega^i) \right] \\
&= \sum_{i=1}^M \Embb_{\omega^i}\left[ S^i(D^i(\overline{\lambda},\omega^i),\omega^i) - \overline{\lambda}^{\Top} D^i(\overline{\lambda},\omega^i) \right].
\end{align*}
The last equality follows from the definition of the demand function $D^i(\pi,\omega^i)$ in \eqref{eq:demand function} as the optimal response of customers to deterministic prices.

The result follows since the tariff $T^*$ induces the same expected total surplus if $\omega \perp \lambda$, \ie
\begin{align*} 
\overline{\sw}(T^*) 
& = \Embb_{\xi} \bigg[ \sum_{i=1}^M S^i(D^{i}(\pi^*,\omega^i),\omega^i) - \lambda^{\Top} D^{i}(\pi^*,\omega^i)\bigg]  \\
& = \sum_{i=1}^M \Embb_{\omega^i} \bigg[ S^i(D^{i}(\overline{\lambda},\omega^i),\omega^i) - \overline{\lambda}^{\Top} D^{i}(\overline{\lambda},\omega^i)\bigg].
\qedadhoc
\end{align*}

\end{document}